\documentclass[11pt,a4]{elsarticle}
\usepackage{latexsym,amsbsy,amssymb,amsmath,amsthm}
\usepackage{enumitem}
\usepackage{epstopdf}
\usepackage{mathrsfs}
\usepackage{color,graphicx}
\usepackage{color}
\usepackage{geometry}
\usepackage{epstopdf}
\usepackage{multirow}
\usepackage[justification=centering]{caption}
\usepackage{subcaption, algorithmicx}
\usepackage[]{algorithm}
\usepackage{longtable}
\usepackage{tabularx}
\usepackage{xfrac}

\usepackage{epstopdf}
\usepackage{graphics}
\usepackage{epsf}
\usepackage{amscd}
\usepackage{wrapfig}
\usepackage[mathscr]{euscript}
\usepackage[english]{babel}
\usepackage{epsfig}
\usepackage{dsfont}
\everymath{\displaystyle}
\numberwithin{equation}{section}
\newtheorem{theorem}{Theorem}[section]
\newtheorem{remark}[theorem]{Remark}

\newcommand{\bv}{\pmb{v}}

\newcommand{\bn}{\pmb{n}}
\newcommand{\bw}{\pmb{w}}

\newcommand{\bx}{\pmb{x}}
\newcommand{\by}{\pmb{y}}

\newcommand{\bD}{\pmb{D}}

\newcommand{\bu}{\pmb{u}}
\newcommand{\btau}{\pmb{\tau}}
\newcommand{\iK}{\mathcal{K}}

\newcommand{\mR}{\mathbb{R}}

\newcommand{\Div}{\text{div}\;}
\newcommand{\iS}{\mathcal{S}}
\newcommand{\iT}{\mathcal{T}}

\newcommand{\iN}{\mathcal{N}}

\newcommand{\iU}{\mathcal{U}}
\newcommand{\iE}{\mathcal{E}}

\newcommand{\bphi}{\pmb{\varphi}}

\begin{document}
\begin{frontmatter}
		
\title{Fully implicit local time-stepping methods for advection-diffusion problems in mixed formulations}
		\author[phuong]{Thi-Thao-Phuong Hoang\fnref{fn}} \ead{tzh0059@auburn.edu}
		
		\address[phuong]{Department of Mathematics and Statistics, Auburn University, Auburn, AL 36849, USA.}
		
		\fntext[fn]{Partially supported by the US National Science Foundation under grant number DMS-1912626.}		
\begin{abstract}
This paper is concerned with numerical solution of transport problems in heterogeneous porous media.  A semi-discrete continuous-in-time formulation of the linear advection-diffusion equation is obtained by using a mixed hybrid finite element method,  in which the flux variable represents both the advective and diffusive flux, and the Lagrange multiplier arising from the hybridization is used for the discretization of the advective term.  Based on global-in-time and nonoverlapping domain decomposition,  we propose two implicit local time-stepping methods to solve the semi-discrete problem.  The first method uses the time-dependent Steklov-Poincar\'e type operator 
and the second uses the optimized Schwarz waveform relaxation (OSWR) with Robin transmission conditions.  For each method, we formulate a space-time interface problem which is solved iteratively. Each iteration involves solving the subdomain problems independently and globally in time; thus,  different time steps can be used in the subdomains.  The convergence of the fully discrete OSWR algorithm with nonmatching time grids is proved.  Numerical results for problems with various Pecl\'et numbers and discontinuous coefficients, including a prototype for the simulation of the underground storage of nuclear waste,  are presented to illustrate the performance of the proposed local time-stepping methods. 
\end{abstract}


\end{frontmatter}

%
%
%
%
\section{Introduction}
\label{Sec:intro}  

Numerical simulations of transport problems in heterogeneous porous media is a subject of great importance in science and engineering.  For applications in hydrology,  various geological layers with different hydrogeological properties are involved in the simulations. Consequently,  the time scales may vary over several order of magnitudes across these layers. This is particularly the case when one simulates the transport of contaminants in and around a nuclear waste repository.  Clearly,  using a single-time step size throughout the entire domain is computationally inefficient; instead,  one should use different time steps in different parts of the domain depending on their physical properties.  In addition,  for the application we consider,  large time step sizes are desirable due to the long time simulation as the nuclear waste decays very slowly.  Therefore,  we propose to use global-in-time,  nonoverlapping domain decomposition (DD) methods in which the dynamic system is decoupled into dynamic subsystems defined on the subdomains (resulting from a spatial decomposition). Then time-dependent problems are solved implicitly in each subdomain at each iteration and the information is exchanged over the space-time interfaces between subdomains.  As a consequence,  different time steps can be used in the subdomains.  
For spatial discretization, we use mixed methods \cite{BF91,RT91} for their mass conservation property and satisfactory performance on heterogeneous problems.  In addition, to handle advection-dominance problems, we employ the (upwind) mixed hybrid finite element method as proposed in~\cite{Radu11, Radu14}, in which the flux variable approximates the total flux (i.e. both diffusive and advective flux).  It was shown that the new mixed hybrid method~\cite{Radu11, Radu14} is fully mass conservative,  as accurate as the standard mixed method~\cite{Dawson09} while it is more efficient in terms of computational cost and robust for problems with high Pecl\'et numbers. 
%

There are basically two types of global-in-time DD methods: the first approach is based on the physical transmission conditions,  for example,  the Dirichlet-Neumann and Neumann-Neumann waveform relaxation methods  \cite{Mandal14, Kwok14, GKM16, GKM20}. The second approach is the Schwarz iteration based on more general transmission conditions such as Robin or Ventcell conditions; an important class of methods in this category is the Optimized Schwarz Waveform Relaxation (OSWR) algorithm  where additional coefficients involved in the transmission conditions are optimized to improve convergence rates \cite{OSWRwave03,Martin05,GH07,GHK07,BlayoHJ,BGH09,HJO10,OSWRDG11,Japhet12,BDF13,HHM13,BGGH16}.  Both approaches were studied with mixed formulations of the pure diffusion problem in~\cite{H13} and the linear advection-diffusion problem in~\cite{H17}.   Particularly in \cite{H13},  a global-in-time preconditioned Schur method (GTP-Schur) and a global-in-time optimized Schwarz method (GTO-Schwarz) were proposed.  Space-time interface problems were derived using, for GTP-Schur,  the time-dependent Dirichlet-to-Neumann (or Steklov-Poincar\'e) operator and,  for GTO-Schwarz,  the time-dependent Robin-to-Robin operator.  

For the GTP-Schur method, the interface problem is solved by preconditioned GMRES with the time-dependent Neumann-Neumann preconditioner,  extended from the Balancing Domain Decomposition (BDD) preconditioner for stationary problems which is known to be efficient when highly heterogeneous coefficients are present \cite{Mandel93,CMW95,MB96}.  Advanced versions of BDD methods are Balancing Domain Decomposition by Constraints (BDDC) methods which introduce a global coarse problem to obtain the condition number bound when the number of subdomains increases (see, e.g., \cite{Cros02, D03, FP03,MDT05,Tu05,LW06,BS07}). 

For the GTO-Schwarz method, the interface problem is solved by either Jacobi iterations or GMRES.  The former choice is equivalent to the OSWR algorithm with Robin transmission conditions, from which the optimized Robin parameters are computed by minimizing the convergence factor in the Fourier transform domain \cite{GH07, GHK07,BGH09,BDF13}.  The interface problems are in mixed form and are space-time for both GTP-Schur and GTP-Schwarz methods, thus nonconforming time grids are possible via a suitable $L^{2}$ projection in time.  An optimal projection algorithm can be found in \cite{OSWRwave03,Projection1d,GJ13}. 

In \cite{H17}, the two methods were extended to the case of advection-diffusion equations with {\em operator splitting} to treat advection and diffusion with different numerical schemes.  For the temporal discretization, the advection is approximated with the explicit Euler method (where sub-time steps are used and constrained by the CFL condition) and the diffusion with the implicit Euler method.  For the spatial discretization, both are approximated with locally mass conservative schemes: the advection with an upwind, cell-center finite volume scheme and the diffusion with a mixed finite element method.  The discrete interface problems for the GTP-Schur and GTO-Schwarz methods are obtained by introducing new unknowns to enforce Dirichlet
transmission conditions between subdomains for the advection step while the diffusion step is handled in the same manner as in \cite{H13}.  In other words, the transmission conditions for the advection part and the diffusion part are separated due to operator splitting.  Consequently,  for the GTO-Schwarz,  the advection plays no role in computing the optimized Robin parameters.  In addition,  it was observed numerically that the GTP-Schur does not perform well when advection is dominant; particularly,  the convergence speed with the (generalized) Neumann-Neumann preconditioner can be even slower than using no preconditioner. 

The objective of this work is to develop fully implicit local time-stepping methods for heterogeneous linear transport problems based on global-in-time DD and the mixed {\em hybrid} finite element method proposed in~\cite{Radu11, Radu14}.  The finite element space is defined using the lowest-order Raviart-Thomas elements in which the total flux is the vector variable,  and the Lagrange multiplier arising in the hybridization is used to discretize the advective term.  Note that for the operator splitting scheme considered in~\cite{H17}, the flux variable approximates the diffusive flux only,  and the advective term is approximated using an upwind operator based on the information from the adjacent elements.  Differently from~\cite{H17}, here we will formulate the fully discrete interface problems for the GTP-Schur and GTO-Schwarz methods where there are no separate interface unknowns for the advection and diffusion.  In addition, unlike \cite{H17} where the advection is treated explicitly and the diffusion implicitly,  the methods proposed in this work are fully implicit with no CFL constraint on the time step size.   
For the GTO-Schwarz method,  the Robin parameters are optimized by taking into account the effects of both advection and diffusion, and we shall prove the convergence of the associated discrete OSWR algorithm with nonconforming time grids.   Note that in this work we focus on the use of local time stepping and only treat conforming spatial discretization.  The reader is referred to \cite{Mortar89, Mortar00, Mortar07, Mortar21}, where mortar mixed methods on nonmatching spatial grids are developed.  
 
The rest of the paper is organized as follows: in the next section we present the model problem and its spatial discretization by the upwind-mixed hybrid finite element method~\cite{Radu11, Radu14}. In Section~\ref{sec:dd}, we formulate two global-in-time decoupling methods using the semi-discrete physical and Robin transmission conditions.  The nonconforming time discretization and the fully discrete interface problems are introduced in Section~\ref{sec:time}; convergence of the OSWR algorithm is also proved where different time steps are used in the subdomains.  In Section~\ref{sec:NumRe}, two-dimensional numerical experiments are carried out to investigate the performance of the proposed methods on different test cases with various Pecl\'et numbers, including one prototype for nuclear waste disposal simulation. 

%
%
%
%
\section{Model problem and its spatial discretization by mixed hybrid finite elements}
\label{sec:model}
For a bounded domain $ \Omega $ of $ \mR^{2}$ with Lipschitz boundary $ \partial \Omega $ and some fixed time $T>0$, consider the following linear advection-diffusion problem
\begin{equation} \label{eq:model}
\begin{array}{rll} \omega \partial_{t} c + \nabla \cdot (\bu c - \bD \nabla c  ) & =f & \text{in} \; \Omega \times (0,T ),\\
c&=0 & \text{on} \; \partial \Omega\times (0,T), \\
c(\cdot, 0) & = c_{0} & \text{in} \; \Omega,
\end{array} 
\end{equation}
where $ c $ is the concentration of a contaminant dissolved in a fluid, $ f $ the source term, $ \omega $ the porosity, $ \bu $ the Darcy velocity ({\em assumed to be given and time independent}), $\bD$ a symmetric time-independent diffusion tensor. For simplicity,  we have imposed only Dirichlet boundary conditions; the analysis presented in the following can be generalized to other types of boundary conditions. Here and throughout the paper, we assume:
\begin{itemize}
\item[(A1)] $\omega$ is bounded above and below by positive constants, $0<\omega_{-} \leq \omega(\bx) \leq \omega_{+}$ for all $\bx \in \Omega$;
\item[(A2)] There exist positive constants $\delta_{-}$ and $\delta_{+}$ such that $\delta_{-} \vert \by\vert^{2} \leq \by^{T} \bD^{-1}(\bx) \by \leq \delta_{+} \vert \by \vert^{2}$ for all $\by \in \mR^{2}$ and $\bx \in \Omega$;
\item[(A3)] $\bu \in (W^{1,\infty}(\Omega))^{2}$, $f \in C(0,T; L^{2}(\Omega))$ and $c \in H^{1}_{0}(\Omega)$.
\end{itemize}

We rewrite \eqref{eq:model} in an equivalent mixed form by introducing the vector field $\bphi$, which consists of both diffusive and advective flux \cite{Radu11, Radu14}:
\begin{equation} \label{eq:mixed}
\begin{array}{rll} 
\omega \partial_{t} c + \nabla \cdot \bphi &=f  & \text{in} \; \Omega \times (0,T), \\
\bphi &= - \bD \nabla c + \bu c & \text{in} \; \Omega \times (0,T),
\end{array} 
\end{equation}
together with the boundary and initial conditions as in \eqref{eq:model}. 

We make use of standard notation for Sobolev spaces and their associated norms  to define the weak formulations and perform convergence analysis in Section~\ref{sec:time}. In particular, we denote by $(\cdot, \cdot)$ the inner product on $L^{2}(\Omega)$ or $(L^{2}(\Omega))^{2}$, and $\| \cdot \|_{l}$ in $H^{l}(\Omega)=W^{l,2}(\Omega)$ (when $l=0$, $H^{l}(\Omega)$ coincides with $L^{2}(\Omega)$. 
For a measurable subset $\Theta \subset \Omega$, we write $(\cdot, \cdot)_{\Theta}$, $\langle \cdot, \cdot \rangle_{\partial \Theta}$ and $\| \cdot \|_{l,\Theta}$ to indicate the inner products or norms considered on $\Theta$. The mixed variational formulation of \eqref{eq:mixed} is given by:

For a.e. $t \in (0,T)$, find $\left (c(t), \bphi(t)\right ) \in L^{2}(\Omega)  \times H(\Div, \Omega)$ such that \vspace{-0.1cm}
\begin{equation} \label{eq:mixedweak}
\begin{array}{rll}
\left (\omega \partial_{t}c, \mu \right ) + \left (\nabla \cdot \bphi, \mu\right ) &= (f,\mu), & \forall \mu \in L^{2}(\Omega), \vspace{3pt}\\
\left (\bD^{-1} \bphi, \bv\right )-\left (\bD^{-1} \bu c, \bv\right ) - \left (\nabla \cdot \bv,c\right ) &=0, & \forall \bv \in H(\text{div}, \Omega).
\end{array} \vspace{-0.2cm}
\end{equation}
Under the assumptions (A1) - (A3), there exists a unique solution $\left (c, \bphi\right ) \in H^{1}\left (0,T; L^{2}(\Omega)\right )  \cap L^{2}(0,T; H^{2}(\Omega) \times L^{2}\left (0,T; H(\Div, \Omega)\right )$ to problem~\eqref{eq:mixedweak} as shown in \cite[Theorem 3.2]{Radu14}.  \vspace{4pt}

%
%
%
%
To find numerical solutions to \eqref{eq:mixedweak}, we use a mixed hybrid finite element (MHFE) method proposed in~\cite{Radu11,Radu14}.  Let $ \iK_{h} $ be a finite element partition of $ \Omega $ into rectangles and let $\iE_h$ be the set of all edges of elements of $ \iK_{h} $: $\iE_{h}:= \iE_{h}^{I} \cup \iE_{h}^{D},$
where $\iE_{h}^{I}$ is the set of all interior edges and $\iE_{h}^{D}$ the set of all edges on the boundary. 
For $K \in \iK_{h} $, let $\bn_K$ be the unit, normal, outward-pointing vector field on the boundary $\partial K$; for each edge $E \subset \partial K$, we denote by $\bn_{E}$ the unit normal vector of $E$, outward to $K$. Let $h_{K}=\text{diam}(K)$ and $h=\max_{K \in \iK_{h}} h_{K}$. 
The MHFE scheme is based on the mixed finite elements together with the hybridization technique, in which the continuity constraint of the normal components of the fluxes over inter-element edges is relaxed via the use of Lagrange multipliers.  The discrete spaces for the scalar and vector variables are defined based on the lowest-order Raviart-Thomas space as
\begin{align*}
M_{h}&:=\left \{ \mu \in L^{2}(\Omega): \mu_{\mid K}=\text{constant}, \; \forall K \in \iK_{h} \right \}, \\
\Sigma_{K} &:= \left \{ \bv:  K \rightarrow \mR^{2}, \; \bv = \left (a_{K} + b_{K} x, a^{\prime}_{K} +b^{\prime}_{K} y \right ), \, (a_{K}, b_{K}, a^{\prime}_{K}, b^{\prime}_{K}) \in \mR^{4} \right \}, \; \text{for} \; K \in \iK_{h}, \\
\Sigma_{h} & := \left \{ \bv \in (L^{2}(\Omega))^{2}: \bv\vert_{K} \in \Sigma_{K} , \; \forall K \in \iK_{h} \right \}  \subset (L^{2}(\Omega))^{2}.
\end{align*}
The discrete space for the Lagrange multiplier representing the trace of the concentration on the edges is given by
$$
\Theta_{h} := \left \{ \theta \in L^{2}(\iE_{h}): \theta\vert_{E}=\, \text{constant on E}, \; \forall E \in \iE_{h} \; \text{and } \, \theta\vert_{E}=0, \; \forall E \in \iE_{h}^{D} \right \}.
$$
For $ c_{h}(t) \in M_{h} $, we have the unique representation $$c_{h} (t,\bx)= \sum_{K \in \iK_{h}} c_{K}(t) \psi_{K}(\bx),$$ where $\psi_{K}$ is the characteristic function of element $K \in \iK_{h}$, and  $c_{K}$ represents the average value of $ c_{h} $ on  $ K $. Similarly, for $\theta_{h} \in \Theta_{h}$, it can be expressed as
$$\theta_{h} (t,\xi)= \sum_{E \in \iE_{h}} \theta_{E}(t) \psi_{E}(\xi),$$ where $\psi_{E}$ is the characteristic function of edge $E \in \iE_{h}$, and  $\theta_{E}$ is the average values of $ \theta_{h} $ on $ E $.
For $\bphi_{h} \in \Sigma_{h}$, the function is defined locally as
$$ \bphi_{h} (t,\bx) \vert_{K}= \sum_{E \subset \partial K} \varphi_{KE} (t)\bw_{KE}(\bx), 
$$
where $\varphi_{KE}$ is the normal flux leaving $K$ through the edge $E$ and $\{\bw_{KE}\}_{E \subset \partial K}$ are the basis functions of the local Raviart-Thomas space $\Sigma_{K}$ satisfying
$$ \int_{E^{\prime}} \bw_{KE} \cdot \bn_{K} = \delta_{E,E^{\prime}}, \; \forall E^{\prime} \subset \partial K. 
$$
We denote by $\bu_{h}$ the projection of $\bu$ on $\Sigma_{h}$ which is defined as
$$\bu_{h}:=\sum_{K \in \iK_{h}} \sum_{E \subset \partial K} u_{KE}\bw_{KE}(\bx), \quad \text{where} \; \, u_{KE} = \frac{1}{\vert E \vert}\int_{E} \bu \cdot \bn_{K}, \; \forall E \subset \partial K, \, \forall K \in \iK_{h}.$$
The mixed hybrid variational formulation for the monodomain problem is given by: 

For a.e. $t \in (0,T)$, find $\left (c_{h}(t), \bphi_{h}(t), \theta_{h}(t)\right ) \in M_{h} \times \Sigma_{h} \times \Theta_{h}$ such that
\begin{equation} \label{eq:dismono}
\begin{array}{rll}
\left (\omega \partial_{t}c_{h}, \mu_{h}\right ) + \left (\nabla \cdot \bphi_{h}, \mu_{h}\right ) &= (f,\mu_{h}), & \forall \mu_{h} \in M_{h}, \vspace{3pt}\\
\left (\bD^{-1} \bphi_{h}, \bv_{h}\right )-\left (\bD^{-1} \bu_{h} c_{h}, \bv_{h}\right ) - \left (\nabla \cdot \bv_{h},c_{h}\right ) + \sum_{K \in \iK_{h}} \left \langle \theta_{h}, \bv_{h} \cdot \bn_{K} \right \rangle_{\partial K} &=0, & \forall \bv_{h} \in \Sigma_{h}, \vspace{3pt}\\
\sum_{K \in \iK_{h}} \left \langle \vartheta_{h}, \bphi_{h} \cdot \bn_{K}\right \rangle_{\partial K} & = 0, & \forall \vartheta \in \Theta_{h}.
\end{array}
\end{equation}
The last equation enforces the continuity of the normal components of the fluxes over inter-element edges so that the vector variable $ \bphi_{h} \in \Sigma_{h}$ belongs to $H(\text{div}, \Omega)$.

For the space-discrete advection term in \eqref{eq:dismono}$_{2}$, instead of using the piecewise constant concentration, we employ the Lagrange multiplier as in \cite{Radu11, Radu14} and obtain the following upwind-mixed scheme:
\begin{equation} \label{eq:upwindmixed}
\begin{array}{rll}
\left (\omega \partial_{t}c_{h}, \mu_{h}\right ) + \left (\nabla \cdot \bphi_{h}, \mu_{h}\right ) &= (f,\mu_{h}), & \forall \mu_{h} \in M_{h}, \vspace{3pt}\\
\left (\bD^{-1} \bphi_{h}, \bv_{h}\right )-\sum_{K \in \iK_{h}} \sum_{E \subset \partial K} u_{KE} \; \iU_{KE}\left (c_{K},\theta_{E}\right ) \left (\bD^{-1} \bw_{KE}, \bv_{h}\right )& &\\
- \left (\nabla \cdot \bv_{h},c_{h}\right ) + \sum_{K \in \iK_{h}} \left \langle \theta_{h}, \bv_{h} \cdot \bn_{K} \right \rangle_{\partial K} &=0, & \forall \bv_{h} \in \Sigma_{h}, \vspace{3pt}\\
\sum_{K \in \iK_{h}} \left \langle \vartheta_{h}, \bphi_{h} \cdot \bn_{K}\right \rangle_{\partial K} & = 0, & \forall \vartheta \in \Theta_{h},
\end{array}
\end{equation}
where the upwind values $\iU_{KE}\left (c_{K},\theta_{E}\right )$ are computed by 
\begin{equation} \label{eq:upwind}
\iU_{KE}\left (c_{K},\theta_{E}\right ) = \left \{ \begin{array}{ll} c_{K}, & \text{if} \; u_{KE} \geq 0, \vspace{3pt}\\
2\theta_{E} - c_{K}, & \text{otherwise}.
\end{array} \right .
\end{equation}
We can also replace the upwind values by the Lagrange multipliers:
\begin{equation} \label{eq:upwind2}
\iU_{KE}\left (c_{K},\theta_{E}\right ) = \theta_{E}.
\end{equation}
Numerical tests \cite{Radu11} suggest that the scheme \eqref{eq:upwindmixed} with the upwind values defined by~\eqref{eq:upwind} is efficient for strongly advection-dominated problems, while using \eqref{eq:upwind2} gives good performance for the case where advection is moderately dominant. 
Both choices preserve the convergence properties of the discretization scheme, and we shall consider formula \eqref{eq:upwind2} in our numerical experiments.
By taking the test functions to be basis functions in \eqref{eq:upwindmixed}, we obtain a system of linear equations for mass conservation, the flux and its continuity over internal edges. 
\paragraph{\bf (i) Mass conservation equation} Let $\omega_{K}$ be the average of $\omega$ on $K$, \eqref{eq:upwindmixed}$_{1}$ implies
\begin{equation} \label{eq:dismass}
\vert K \vert \omega_{K} \partial_{t} c_{K} + \sum_{E \subset \partial K} \varphi_{KE} =\int_{K} f \, d\bx, \quad \forall K \in \iK_{h}.
\end{equation}
\paragraph{\bf (ii) Equation for the flux} Denote by $A_{KEE^{\prime}} = \int_{K} (\bD^{-1} \bw_{KE^{\prime}}) \cdot \bw_{KE} \, d\bx$, for $E, E^{\prime} \subset \partial K$, $K \in \iK_{h}$, \eqref{eq:upwindmixed}$_{2}$ becomes
\begin{equation} \label{eq:disflux}
\sum_{E^{\prime} \subset \partial K} A_{KEE^{\prime}} \varphi_{KE^{\prime}} - \sum_{E^{\prime} \subset \partial K} A_{KEE^{\prime}} u_{KE^{\prime}} \; \iU_{KE^{\prime}}\left (c_{K},\theta_{E^{\prime}}\right ) -c_{K} + \theta_{E}=0, \quad \forall K \in \iK_{h}, \forall E \subset \partial K.
\end{equation}
\paragraph{\bf (iii) Continuity of the flux over internal edges} It is deduced from \eqref{eq:upwindmixed}$_{3}$ that
\begin{equation} \label{eq:disLag}
\sum_{\substack{K \in \iK_{h} \\ E \subset \partial K}} \varphi_{KE} = 0, \quad \forall E \in \iE_{h}^{I}.
\end{equation}
The system \eqref{eq:dismass}-\eqref{eq:disLag} is completed with some appropriate boundary conditions and initial condition. 
\begin{remark}
The fully discrete problem obtained by discretizing  \eqref{eq:upwindmixed} in time using the implicit Euler method was analyzed in \cite[Theorem 4.4]{Radu14}.  Expected error estimates were proved, given that the choice of upwind values satisfies the following inequality 
$$ \vert \iU_{KE}\left (c_{K},\lambda_{E}\right ) - c_{K} \vert \leq C \vert \lambda_{E} - c_{K} \vert, \quad \forall E \subset \partial K, \, \forall K \in \iK_{h},
$$
for some constant $C$ independent of the mesh size and time step size. Note that the upwind values given by either \eqref{eq:upwind} or \eqref{eq:upwind2} satisfy this requirement.
\end{remark}
%
%
%
\section{Semi-discrete, global-in-time domain decomposition methods} \label{sec:dd}
We use nonoverlapping domain decomposition methods that allow different time steps to solve the semi-discrete problem~\eqref{eq:dismono}.
For simplicity, we consider a decomposition of $ \Omega $ into two nonoverlapping subdomains
$ \Omega_{1} $ and $ \Omega_{2} $ separated by an interface $\Gamma$:
$$
\Omega_{1} \cap \Omega_{2} = \emptyset;\quad  \Gamma
= \partial \Omega_{1} \cap \partial \Omega_{2} \cap \Omega, \quad \Omega=  \Omega_{1} \cup \Omega_{2} \cup\Gamma.
$$
The analysis can be generalized to the case of multiple subdomains (see Section~\ref{sec:NumRe}). We assume that the partitions $ \iK_{h,1} $ of subdomain $ \Omega_{1} $ and $ \iK_{h,2} $ of subdomain $ \Omega_{2} $ are such that their
union $ \iK_{h}= \bigcup_{i=1}^2 \iK_{h,i}$ forms a finite element partition of $\Omega $. Let $\iE_{h,i}$ be the set of all edges of elements of $ \iK_{h,i}$ and $\iE_{h,i}^{D}$ the set of edges on the external boundary $\partial \Omega_{i} \cap \partial \Omega$, for $i=1,2$. Notice that in this work we focus on the use of different time steps, and assume the spatial discretization is conforming.  Thus we denote by $ \iE_{h}^{\Gamma} $ the set of edges of elements of $\iK_{h,1}$ or $\iK_{h,2}$ that lie on $ \Gamma$.
For $i=1,2$, let $ \bn_{i} $ denote the unit, normal, outward-pointing vector field on $\partial\Omega_i$, and for any scalar, vector or tensor-valued function $ \psi $ defined on $\Omega$, let $\psi_i$ denote the restriction  of $ \psi$ to $ \Omega_{i} $. 

Let $ M_{h} $, $ \Sigma_{h} $ and $\Theta_{h}$ denote the mixed finite element spaces and the Lagrange multiplier space as defined in Section~\ref{sec:model}, and let $ M_{h,i} $, $ \Sigma_{h,i} $ and $\Theta_{h,i}$, $ i=1,2, $ be the spaces of restrictions of the functions in these spaces to $ \Omega_{i} $.  In particular:
$$ \Theta_{h,i} := \left \{ \theta \in L^{2}(\iE_{h,i}): \theta\vert_{E}=\, \text{constant on E}, \; \forall E \in \iE_{h,i} \; \text{and } \, \theta\vert_{E}=0, \; \forall E \in \iE_{h,i}^{D}\right \}. 
$$
In addition, we define the space
$$ \Theta_{h,i}^{\Gamma,0} := \left \{ \theta \in  \Theta_{h,i}:  \theta\vert_{E}=0, \; \forall E \in \iE_{h}^{\Gamma}\right \}. 
$$
to take into account the interface as part of the subdomain boundary.
As for the single domain case, we have the following representations of functions in the subdomain finite element spaces:
\begin{equation}
\begin{array}{ll}
c_{h,i} (t,\bx)&= \sum_{K \in \iK_{h,i}} c_{i,K}(t) \psi_{K}(\bx), \quad \theta_{h,i} (t,\xi)= \sum_{E \in \iE_{h,i}} \theta_{i,E}(t) \psi_{E}(\xi), \\
 \bphi_{h,i} (t,\bx) &= \sum_{K \in \iK_{h,i}} \sum_{E \subset \partial K} \varphi_{i,KE} (t)\bw_{KE}(\bx).
\end{array}
\end{equation}
Consequently, we may identify functions $c_{h,i} \in M_{h,i}$ with vectors $\left ( c_{i,K}\right )_{K \in \iK_{h,i}}$, $\theta_{h,i} \in \Theta_{h,i}$ with $\left (\theta_{i,E}\right )_{E \in \iE_{h,i}}$, and $\bphi_{h,i} \in \Sigma_{h,i}$ with $\left (\varphi_{i,KE}\right )_{K \in \iK_{h,i}, E \subset K} $ when necessary to simplify the presentation. 
To define the transmission conditions between the subdomains, we introduce the following interface space consisting of piecewise constant functions on $\Gamma$:
\begin{equation}
\Lambda_{h} = \left \{ \lambda \in L^{2}(\iE_{h}^{\Gamma}): \lambda\vert_{E}=\, \text{constant on E}, \; \forall E \in \iE_{h}^{\Gamma}  \right \}.
\end{equation}


With the notation introduced above,  the monodomain problem~\eqref{eq:dismono} is equivalent to the following time-dependent subdomain problems in $\Omega_{1}$ and $\Omega_{2}$:
\begin{equation} \label{eq:multi}
\begin{array}{ll}
&\left (\omega_{i} \partial_{t}c_{h,i}, \mu_{h}\right )_{\Omega_{i}} + \left (\nabla \cdot \bphi_{h,i}, \mu_{h}\right )_{\Omega_{i}} = (f,\mu_{h})_{\Omega_{i}}, \quad \forall \mu_{h} \in M_{h,i}, \vspace{3pt}\\
&\left (\bD_{i}^{-1} \bphi_{h,i}, \bv_{h}\right )_{\Omega_{i}}-\sum_{K \in \iK_{h,i}} \sum_{E \subset \partial K} u_{i,KE} \; \iU_{KE}\left (c_{i,K},\theta_{i,E}\right ) \left (\bD_{i}^{-1} \bw_{KE}, \bv_{h}\right )_{\Omega_{i}}  \\
& \hspace{3cm} -\left (\nabla \cdot \bv_{h},c_{h,i}\right )_{\Omega_{i}} + \sum_{K \in \iK_{h,i}} \left \langle \theta_{h,i}, \bv_{h} \cdot \bn_{K} \right \rangle_{\partial K} =0, \; \forall \bv_{h} \in \Sigma_{h,i}, \vspace{3pt}\\
&\sum_{K \in \iK_{h}} \left \langle \vartheta_{h}, \bphi_{h,i} \cdot \bn_{K}\right \rangle_{\partial K} = 0, \quad \forall \vartheta_{h} \in \Theta_{h,i}^{\Gamma, 0},
\end{array}
\end{equation}
for $i=1,2$, together with the time-dependent transmission conditions:
\begin{equation} \label{eq:tcs}
\begin{array}{c}
 \left \langle \theta_{h,1}(t), \eta_{h}  \right \rangle_{\Gamma}  = \left \langle \theta_{h,2}(t), \eta_{h}  \right \rangle_{\Gamma}, \\
\left \langle \eta_{h} , \bphi_{h,1}(t) \cdot \bn_{1} + \bphi_{h,2}(t) \cdot \bn_{2} \right \rangle_{\Gamma} = 0,
\end{array} \quad \forall \eta_{h} \in \Lambda_{h},
\end{equation}
for $t \in (0,T)$, where $\langle \cdot, \cdot \rangle_{\Gamma}$ denotes the inner product of $L^{2}(\Gamma)$.
Alternatively, one may use the following equivalent Robin transmission conditions:
\begin{equation} \label{eq:robintcs}
\left. \begin{array}{ll} \left \langle -\bphi_{h,1}(t) \cdot \bn_{1} + \alpha_{1,2}  \theta_{h,1}(t), \eta_{h}  \right \rangle_{\Gamma} & = \left \langle -\bphi_{h,2}(t) \cdot \bn_{1} + \alpha_{1,2}  \theta_{h,2}(t), \eta_{h} \right \rangle_{\Gamma},\\
\left \langle -\bphi_{h,2}(t) \cdot \bn_{2} + \alpha_{2,1}  \theta_{h,2}(t), \eta_{h} \right \rangle_{\Gamma}& = \left \langle -\bphi_{h,1}(t) \cdot \bn_{2} + \alpha_{2,1}  \theta_{h,1}(t), \eta_{h} \right \rangle_{\Gamma},
\end{array} \right . \forall \eta_{h} \in \Lambda_{h},
\end{equation}
where $\alpha_{1,2}$ and $\alpha_{2,1}$ are a pair of positive parameters. 
By the definition of the finite dimensional approximation spaces, the transmission condition in \eqref{eq:tcs} can be rewritten as
%
\begin{equation} \label{eq:PhysTCs}
\begin{array}{c}
 \theta_{1,E} (t)= \theta_{2,E}(t),  \\
\varphi_{1,KE}(t) + \varphi_{2, K^{\prime}E}(t) = 0, 
 \end{array} \quad  \forall \; E \in \iE_{h}^{\Gamma}, \; t \in (0,T), 
\end{equation}
where $E = \partial K \cap \partial K^{\prime}, \; \text{with} \; K \in \iK_{h,1} \; \text{and} \; K^{\prime} \in \iK_{h,2}$ (i.e. $K$ and $K^{\prime}$ are elements in different subdomains sharing the interface edge $E$). 
Similarly, the Robin transmission conditions \eqref{eq:robintcs} are equivalent to
\begin{equation} \label{eq:RobinTCs}
\left. \begin{array}{ll} -\varphi_{1,KE}(t) + \alpha_{1,2}  \theta_{1,E} (t)& = \varphi_{2, K^{\prime}E}(t)+ \alpha_{1,2} \theta_{2,E}(t)\\
- \varphi_{2, K^{\prime}E}(t) + \alpha_{2,1} \theta_{2,E}(t) & = \varphi_{1,KE}(t) + \alpha_{2,1} \theta_{1,E}(t)
\end{array} \right . \quad \forall E \in \iE_{h}^{\Gamma}, \; t \in (0,T).
\end{equation}
The first method that we consider is based on \eqref{eq:multi} together with the physical transmission conditions \eqref{eq:PhysTCs} while the second method is based on \eqref{eq:multi} together with the Robin transmission conditions \eqref{eq:RobinTCs}.   For the latter method, the parameters $\alpha_{i,j}$ may be optimized to improve
the convergence rate of the iterative scheme.  For precise details of how the optimization is carried out, see \cite{GHK07,Japhet12}. 

\subsection{Method 1: Global-in-Time Preconditioned Schur (GTP-Schur)}
We impose the first transmission condition in \eqref{eq:tcs} as Dirichlet boundary conditions on the interface of the subdomain problems:
\begin{equation}  \label{eq:DirBcIf}
  \theta_{h,i}(t)\vert_{\Gamma} = \lambda_{h}(t), \quad t \in (0,T), \; \text{for} \; i=1,2,
\end{equation}
for some given $\lambda_{h}(t) \in \Lambda_{h}$. We introduce the space
\begin{align}
\Theta_{h,i}^{\Gamma, \lambda_{h}} &:= \big \{ \theta \in L^{2}(\iE_{h,i}): \theta\vert_{E}=\, \text{constant on E}, \; \forall E \in \iE_{h,i}, \theta\vert_{E}=0, \; \forall E \in \iE_{h,i}^{D}, \nonumber\\
& \hspace{8cm} \text{and } \, \theta\vert_{E}=\lambda_{h} \vert_{E}, \; \forall E \in \iE_{h,i}^{\Gamma}\big \}, 
\end{align}
for the Lagrange multipliers with nonhomogeneous Dirichlet conditions on $\Gamma$. We denote by $$\left (c_{h,i}(\lambda_{h},f,c_{0}), \bphi_{h,i}(\lambda_{h},f,c_{0}), \theta_{h,i}(\lambda_{h},f,c_{0})\right ) \in  H^{1}(0,T;M_{h,i}) \times L^{2}(0,T;\Sigma_{h,i}) \times L^{2}(0,T;\Theta_{h,i}^{\Gamma, \lambda_{h}})$$ the solution to the time-dependent subdomain problem with Dirichlet interface condition~\eqref{eq:DirBcIf}:
\begin{equation} \label{eq:disM1}
\begin{array}{ll}
&\left (\omega_{i} \partial_{t}c_{h,i}, \mu_{h}\right )_{\Omega_{i}} + \left (\nabla \cdot \bphi_{h,i}, \mu_{h}\right )_{\Omega_{i}} = (f,\mu_{h})_{\Omega_{i}}, \quad \forall \mu_{h} \in M_{h,i}, \vspace{3pt}\\
&\left (\bD_{i}^{-1} \bphi_{h,i}, \bv_{h}\right )_{\Omega_{i}}-\sum_{K \in \iK_{h,i}} \sum_{E \subset \partial K} u_{i,KE} \; \iU_{KE}\left (c_{i,K},\theta_{i,E}\right ) \left (\bD_{i}^{-1} \bw_{KE}, \bv_{h}\right )_{\Omega_{i}}  \\
& \hspace{2cm } - \left (\nabla \cdot \bv_{h},c_{h,i}\right )_{\Omega_{i}}+ \sum_{K \in \iK_{h,i}} \left \langle \theta_{h,i}, \bv_{h} \cdot \bn_{K} \right \rangle_{\partial K} =0, \quad  \forall \bv_{h} \in \Sigma_{h,i}, \vspace{3pt}\\
& \sum_{K \in \iK_{h}} \left \langle \vartheta_{h}, \bphi_{h,i} \cdot \bn_{K}\right \rangle_{\partial K} = 0, \quad \forall \vartheta_{h} \in \Theta_{h,i}^{\Gamma, 0}, \vspace{3pt}\\
& (c_{h,i}(0), \mu_{h})_{\Omega_{i}} = (c_{0}, \mu_{h})_{\Omega_{i}}, \quad \forall \mu_{h} \in M_{h,i},
\end{array}
\end{equation}
for $i=1,2$.  

The semi-discrete interface problem for Method 1 is obtained by enforcing the remaining transmission condition, i.e. the flux continuity equation. Toward that end, we define the interface operators, for $i=1,2$:
\begin{equation*}
\begin{array}{lcll}
	\iS_{h,i}\,:& L^{2}(0,T; \Lambda_{h}) & \longrightarrow & L^{2}(0,T;\Lambda_{h}), \\
	& \lambda_{h} & \mapsto &  -\bphi_{h,i}(\lambda_{h},0,0) \cdot \bn_{i}\vert_{\Gamma}, \vspace{0.2cm}\\
\end{array}
\end{equation*}
and
\begin{equation*}
\begin{array}{lcll}
	\chi_{h,i}\,:& L^{2}(0,T; \Omega_{i}) \times H^{1}_{\ast}(\Omega_{i})& \longrightarrow & L^{2}(0,T;\Lambda_{h}), \\
	& (f,c_{0}) & \mapsto &  \bphi_{h,i}(0,f,c_{0}) \cdot \bn_{i}\vert_{\Gamma}, \vspace{0.2cm}\\
\end{array}
\end{equation*}
where $H^{1}_{\ast}(\Omega_{i})=\{ \mu \in H^{1}(\Omega_{i}): v=0 \; \text{on} \; \partial \Omega_{i} \cap \partial \Omega\}$. The operators $\iS_{h,i}$, $i=1,2,$ are the (time-dependent) Schur operators. Letting $\iS_{h}=\iS_{h,1}+\iS_{h,2}$ and $\chi_{h}=\chi_{h,1}+\chi_{h,2}$, we rewrite problem~\eqref{eq:multi}-\eqref{eq:tcs} as a space-time interface problem:

Find $\lambda_{h} \in L^{2}(0,T; \Lambda_{h}),$ such that
\begin{equation}\label{eq:IP_M1}
\int_{0}^{T} \int_{\Gamma} \left (\iS_{h}\lambda_{h} \right ) \, \eta_{h}  \, d\gamma \, dt=\int_{0}^{T} \int_{\Gamma} \chi_{h} \, \eta_{h}  \, d\gamma \, dt, \quad \forall \eta_{h}  \in L^{2}(0,T;\Lambda_{h}).
\end{equation}

The interface problem \eqref{eq:IP_M1} is solved iteratively,  and to enhance the convergence of the iteration, we apply a Neumann-Neumann type preconditioner similarly to the one proposed for the pure diffusion case \cite{H13}.  For a function $\phi_{h} \in L^{2}(0,T;\Lambda_{h})$, we define the following spaces to handle Neumann boundary conditions involved in the preconditioner:
\begin{align*}
\Sigma_{h,i}^{\Gamma,\phi_{h}} &= \left \{ \bv \in \Sigma_{h,i}: v_{KE} = \phi_{h} \vert_{E}, \quad \forall E \in \iE_{h}^{\Gamma}, \, K \in \iK_{h,i} \; \text{such that} \, E \subset \partial K \right \}, \;  \\
\Sigma_{h,i}^{\Gamma,0} &= \left \{ \bv \in \Sigma_{h,i}: v_{KE}= 0, \quad \forall E \in \iE_{h}^{\Gamma}, \, K \in \iK_{h,i} \; \text{such that} \, E \subset \partial K \right \}.
\end{align*}
The Neumann-Neumann preconditioner for \eqref{eq:IP_M1} is given by:
\begin{equation} \label{eq:precondS}
(\sigma_{1,2} \iN_{h,1} + \sigma_{2,1} \iN_{h,2}) \iS_{h} \lambda_{h} = \widetilde{\chi_{h}}, 
\end{equation}
where $\sigma_{i,j}$ are some weights such that $\sigma_{1,2} + \sigma_{2,1}=1$, and $\iN_{h,i}$ is a (pseudo-)inverse operator of $\iS_{h,i}$ defined, for $\phi_{h} \in L^{2}(0,T;\Lambda_{h})$, as
\begin{equation} \label{eq:Sinv}
\iN_{h,i} \psi_{h} = \theta_{h,i}(\phi_{h})\vert_{\Gamma}, \quad i=1,2,
\end{equation}
where $\left (c_{h,i}(\phi_{h}), \bphi_{h,i}(\phi_{h}), \theta_{h,i}(\phi_{h})\right ) \in H^{1}(0,T;M_{h,i}) \times L^{2}(0,T;\Sigma_{h,i}^{\Gamma,\phi_{h}}) \times L^{2}(0,T;\Theta_{h,i}) $, $i=1,2$, is the solution to the time-dependent subdomain problem with Neumann boundary conditions on the interface and zero initial condition:
\begin{equation} \label{eq:precondsub}
\begin{array}{ll}
&\left (\omega_{i} \partial_{t}c_{h,i}, \mu_{h}\right )_{\Omega_{i}} + \left (\nabla \cdot \bphi_{h,i}, \mu_{h}\right )_{\Omega_{i}} = 0, \quad \forall \mu_{h} \in M_{h,i}, \vspace{3pt}\\
&\left (\bD_{i}^{-1} \bphi_{h,i}, \bv_{h}\right )_{\Omega_{i}}-\sum_{K \in \iK_{h,i}} \sum_{E \subset \partial K} u_{i,KE} \; \iU_{KE}\left (c_{i,K},\theta_{i,E}\right ) \left (\bD_{i}^{-1} \bw_{KE}, \bv_{h}\right )_{\Omega_{i}}  \\
& \hspace{2cm}- \left (\nabla \cdot \bv_{h},c_{h,i}\right )_{\Omega_{i}} + \sum_{K \in \iK_{h,i}} \left \langle \theta_{h,i}, \bv_{h} \cdot \bn_{K} \right \rangle_{\partial K} =0, \quad \forall \bv_{h} \in \Sigma_{h,i}^{\Gamma,0}, \vspace{3pt}\\
& \sum_{K \in \iK_{h}} \left \langle \vartheta_{h}, \bphi_{h} \cdot \bn_{K}\right \rangle_{\partial K} = 0, \quad \forall \vartheta_{h} \in \Theta_{h,i}^{\Gamma, 0}.
\end{array}
\end{equation}
The performance of the preconditioner will be investigated numerically in Section~\ref{sec:NumRe}.  
%
%
%
%
\subsection{Method 2: Global-in-Time Optimized Schwarz (GTO-Schwarz)}
We impose Robin transmission conditions~\eqref{eq:RobinTCs} as boundary conditions on $\Gamma \times (0,T)$ for solving the subdomain problems:
\begin{equation} \label{eq:RobinBcsIf}
\left (-\bphi_{h,i}(t) \cdot \bn_{i} + \alpha_{1,2}  \theta_{h,i}(t)\right )\vert_{\Gamma} = \zeta_{h,i}(t), \quad t \in (0,T),
\end{equation} 
for given $\zeta_{h,i} \in L^{2}(0,T;\Lambda_{h})$ and for $i=1,2$. Denote by $$\left (c_{h,i}(\zeta_{h,i},f,c_{0}), \bphi_{h,i}(\zeta_{h,i},f,c_{0}), \theta_{h,i}(\zeta_{h,i},f,c_{0})\right ) \in H^{1}(0,T;M_{h,i}) \times L^{2}(0,T;\Sigma_{h,i}) \times L^{2}(0,T;\Theta_{h,i})$$ the solution to the time-dependent subdomain problem with Robin interface condition~\eqref{eq:RobinBcsIf}: \vspace{-0.2cm}
\begin{equation} \label{eq:disM2}
\begin{array}{l}
\left (\omega_{i} \partial_{t}c_{h,i}, \mu_{h}\right )_{\Omega_{i}}  + \left (\nabla \cdot \bphi_{h,i}, \mu_{h}\right )_{\Omega_{i}}  = (f,\mu_{h})_{\Omega_{i}} , \; \forall \mu_{h} \in M_{h,i}, \vspace{3pt}\\
\left (\bD_{i}^{-1} \bphi_{h,i}, \bv_{h}\right )_{\Omega_{i}} -\sum_{K \in \iK_{h,i}} \sum_{E \subset \partial K} u_{i,KE} \; \iU_{KE}\left (c_{i,K},\theta_{i,E}\right ) \left (\bD_{i}^{-1} \bw_{KE}, \bv_{h}\right )_{\Omega_{i}}  \\
\hspace{2cm} - \left (\nabla \cdot \bv_{h},c_{h,i}\right )_{\Omega_{i}}  + \sum_{K \in \iK_{h}} \left \langle \theta_{h,i}, \bv_{h} \cdot \bn_{K} \right \rangle_{\partial K} = 0, \quad \forall \bv_{h} \in \Sigma_{h,i}, \vspace{3pt}\\
\sum_{K \in \iK_{h}} \left \langle \vartheta_{h}, \bphi_{h,i} \cdot \bn_{K}\right \rangle_{\partial K}  = 0, \; \forall \vartheta_{h} \in \Theta_{h}^{\Gamma, 0}, \vspace{3pt}\\
 \left \langle -\bphi_{h,i} \cdot \bn_{i}+\alpha_{i,j}  \theta_{h,i}, \eta_{h} \right \rangle_{\Gamma} = \left \langle \zeta_{h,i}, \eta_{h}  \right \rangle_{\Gamma}, \quad \forall \eta_{h} \in \Lambda_{h}, \vspace{3pt}\\
  (c_{h,i}(0), \mu_{h})_{\Omega_{i}} = (c_{0}, \mu_{h})_{\Omega_{i}}, \quad \forall \mu_{h} \in M_{h,i}. 
\end{array}
\end{equation}

We define the interface operators:
\begin{equation*}
\begin{array}{lcll}
	\iS^{R}_{h}\,:& \left (L^{2}(0,T; \Lambda_{h})\right )^{2} &\longrightarrow& \left (L^{2}(0,T; \Lambda_{h})\right )^{2}, \\
	& \left (\begin{array}{c} \zeta_{h,1} \\ \zeta_{h,2} \end{array}\right )	 & \mapsto & \hspace{-0.3cm} \left (\begin{array}{c} \zeta_{h,1} - \left (-\bphi_{h,2}(\zeta_{h,2},0,0) \cdot \bn_{1} + \alpha_{1,2}  \theta_{h,2}(\zeta_{h,2},0,0)\right )\vert_{\Gamma} \\ 
	\zeta_{h,2}-\left (-\bphi_{h,1}(\zeta_{h,1},0,0) \cdot \bn_{2} + \alpha_{2,1}  \theta_{h,1}(\zeta_{h,1},0,0)\right )\vert_{\Gamma}
	 \end{array} \right ),
\end{array} \vspace{-0.2cm}
\end{equation*}
and 
\begin{equation*}
\begin{array}{lcll}
	\chi^{R}_{h}\,:& L^{2}(0,T; \Omega_{i}) \times H^{1}_{\ast}(\Omega_{i})& \longrightarrow & \left (L^{2}(0,T; \Lambda_{h})\right )^{2}, \\
	& (f,c_{0}) & \mapsto &  \hspace{-0.3cm} \left (\begin{array}{c}  \left (-\bphi_{h,2}(0,f,c_{0}) \cdot \bn_{1} + \alpha_{1,2}  \theta_{h,2}(0,f,c_{0})\right )\vert_{\Gamma} \\ 
	\left (-\bphi_{h,1}(0,f,c_{0}) \cdot \bn_{2} + \alpha_{2,1}  \theta_{h,1}(0,f,c_{0})\right )\vert_{\Gamma}
	 \end{array} \right ),
\end{array}
\end{equation*}
for $i=1,2$. 
The semi-discrete interface problem for Method 2 is given by: 

Find $(\zeta_{h,1}, \zeta_{h,2}) \in  \left (H^{1}(0,T; \Lambda_{h})\right )^{2}$ such that
\begin{equation}\label{eq:IP_M2}
\int_{0}^{T}  \int_{\Gamma} S^{R}_{h}\left (\begin{array}{c} \zeta_{h,1} \\ \zeta_{h,2} \end{array}\right ) \cdot \left (\begin{array}{c} \eta_{h,1}  \\ \eta_{h,2}  \end{array}\right )\, d\gamma \, dt=\int_{0}^{T}  \int_{\Gamma} \chi^{R}_{h}\cdot \left (\begin{array}{c} \eta_{h,1}  \\ \eta_{h,2} \end{array}\right ) \, d\gamma \, dt, \; \forall (\eta_{h,1}, \eta_{h,2}) \in \left (L^{2}(0,T;\Lambda_{h})\right )^{2}.
\end{equation}
As for GTP-Schur, we solve the interface problem \eqref{eq:IP_M2} iteratively using Jacobi iterations or GMRES.  The former choice is equivalent to the following {\em Optimized Schwarz waveform relaxation (OSWR) algorithm}: at the $k$th iteration, solve in parallel the subdomain problems:
\begin{equation} \label{eq:semiOSWR}
\begin{array}{l}
\left (\omega_{i} \partial_{t}c_{h,i}^{k}, \mu_{h}\right )_{\Omega_{i}}  + \left (\nabla \cdot \bphi_{h,i}^{k}, \mu_{h}\right )_{\Omega_{i}}  = (f,\mu_{h})_{\Omega_{i}} , \; \forall \mu_{h} \in M_{h,i}, \vspace{3pt}\\
\left (\bD_{i}^{-1} \bphi_{h,i}^{k}, \bv_{h}\right )_{\Omega_{i}} -\sum_{K \in \iK_{h,i}} \sum_{E \subset \partial K} u_{i,KE} \; \iU_{KE}\left (c_{i,K}^{k},\theta_{i,E}^{k}\right ) \left (\bD_{i}^{-1} \bw_{KE}, \bv_{h}\right )_{\Omega_{i}}  \\
\hspace{2cm}  - \left (\nabla \cdot \bv_{h},c_{h,i}^{k}\right )_{\Omega_{i}} + \sum_{K \in \iK_{h}} \left \langle \theta_{h,i}^{k}, \bv_{h} \cdot \bn_{K} \right \rangle_{\partial K}=0, \quad \forall \bv_{h} \in \Sigma_{h,i}, \vspace{3pt}\\
\sum_{K \in \iK_{h}} \left \langle \vartheta_{h}, \bphi_{h,i}^{k} \cdot \bn_{K}\right \rangle_{\partial K}  = 0, \; \forall \vartheta_{h} \in \Theta_{h}^{\Gamma, 0}, \vspace{3pt}\\
 \left \langle -\bphi_{h,i}^{k} \cdot \bn_{i}+\alpha_{i,j}  \theta_{h,i}^{k}, \eta_{h} \right \rangle_{\Gamma} = \left \langle-\bphi_{h,j}^{k-1} \cdot \bn_{i}+\alpha_{i,j}  \theta_{h,j}^{k-1}, \eta_{h}  \right \rangle_{\Gamma}, \quad \forall \eta_{h} \in \Lambda_{h},
\end{array}
\end{equation}
for $i=1,2$,
%
with given initial guesses: $g_{i,j}(t):=-\varphi_{h,j}^{0} \cdot \bn_{i} + \alpha_{i,j}  \theta_{h,j}^{0} \in \Lambda_{h}$ for $i=1,2,$ $j=3-i$ to start the first iterate.  We will prove the convergence of the fully discrete OSWR algorithm with different subdomain time steps in the next section.  
%
%
%
%
\section{Nonconforming time discretizations}
\label{sec:time}
\begin{figure}[h]
\vspace{1.2cm}
\centering
\begin{minipage}[c]{0.5 \linewidth}
\setlength{\unitlength}{1pt} 
\begin{picture}(140,70)(0,0)
\thicklines
\put(0,3){\includegraphics[scale=0.55]{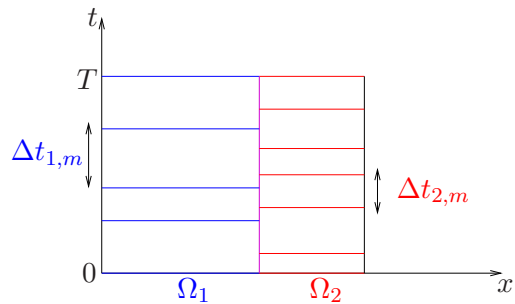} \\}
\put(-1,1){$ 0 $}
\put(-3,74){$ T $}
\put(35,-5){\textcolor{blue}{$ \Omega_{1} $}}
\put(85,-5){\textcolor{red}{$ \Omega_{2}$}}
\put(-28,47){\textcolor{blue}{$\Delta t_{1,m} $}}
\put(118,33){\textcolor{red}{$ \Delta t_{2,m} $}}
\put(156,-3){$ x $}
\put(1,98){$ t $}
\end{picture}
\end{minipage}
\caption{Nonconforming time grids in the subdomains.}
\label{Fig:Time} 
\end{figure}

The GTP-Schur and GTO-Schwarz methods involve solving the subdomain problems globally in time, thus independent time discretizations can be used in the subdomains.  We consider fully discrete problems with nonconforming time grids. Let $ \iT_{1} $ and $ \iT_{2} $
be two possibly different partitions of the time interval $ (0,T) $ into sub-intervals (see Figure~\ref{Fig:Time}).
We denote by $ J_{i,m} $ the time interval $ (t_{i,m}, t_{i,m-1}] $ and by
$ \Delta t_{i,m} := (t_{i,m} - t_{i,m-1}) $ for $ m=1, \hdots, M_{i} $ and $ i=1,2 $, where for simplicity of exposition we have again supposed that we have only two subdomains.
We use the the backward Euler method for discretizing in time, the same idea can be generalized to higher order methods. 
\noindent We denote by $ P_{0}(\mathcal{T}_{i}, W) $ the space of piecewise constant functions in time on grid
$ \mathcal{T}_{i} $ with values in $ W $:
\begin{equation} \label{P0space}
P_{0}(\mathcal{T}_{i}, W) = \left \{ \psi: (0,T) \rightarrow W,\; 
\psi \; \text{ is constant on} \;  J_{i,m}, \ \forall m=1, \dots, M_{i} \right \}. \vspace{-0.1cm}
\end{equation}
In order to exchange data on the space-time interface between different time grids, we define the following
$ L^{2} $ projection $ \Pi_{ji} $ from  $ P_{0} (\mathcal{T}_{i}, W) $ onto $ P_{0}(\mathcal{T}_{j},W) $
(see \cite{OSWRwave03,Japhet12})~: for $ \psi \in P_{0} (\mathcal{T}_{i}, W)$,
$ \Pi_{ji} \psi \hspace{-2pt} \mid_{J_{j,m}} $ is the average value of $ \psi $ on $ J_{j,m} $,
for $ m=1, \dots, M_{j} $: 
\begin{equation} \label{ProjectionTime}
\Pi_{ji} \left ( \psi \right )\mid_{J_{j,m}}
  = \frac{1}{\mid J_{j,m}\mid} \sum_{l=1}^{M_{i}} \int_{J_{j,m} \cap J_{i,l}} \psi. 
\end{equation}
We use the algorithm described in \cite{Projection1d} for effectively performing this projection.

\subsection{For Method 1 (GTP-Schur):} As there is only one unknown $\lambda_{h}$ on the interface, we need
to choose $\lambda_{h}$ piecewise constant in time on one grid, either $\iT_1$ or $\iT_2$. For instance, let $\lambda_{h} = (\lambda_{h}^{m})_{m=1}^{M_{1}} \in P_{0}(\mathcal{T}_{1}, \Lambda_{h})$. The fully discrete counterpart of \eqref{eq:IP_M1} is weakly enforced over the time intervals of $\iT_{1}$ as follows: $\forall m=1, \hdots, M_{1},$
\begin{equation} \label{eq:fullIP_M1}
\begin{array}{c}
\hspace{-2.5cm} \int_{J_{1,m}} \int_{\Gamma}  \left (-\bphi_{h,1}(\lambda_{h},0,0) \cdot \bn_{1}-\Pi_{12}\left (\bphi_{h,2}(\Pi_{21}(\lambda_{h}),0,0)\right ) \cdot \bn_{1}\right ) \eta_{h} \, d\gamma \, dt= \vspace{0.2cm} \\
 \int_{J_{1,m}} \int_{\Gamma}  \left (\bphi_{h,1}(0,f,c_{0}) \cdot \bn_{1}-\Pi_{12}\left (\bphi_{h,2}(0,f,c_{0}) \cdot \bn_{1}\right )\right ) \eta_{h}\, d\gamma \, dt, \quad \forall \eta_{h} \in \Lambda_{h},
 \end{array}
\end{equation}
where $\bphi_{h,i}(\lambda_{h},c_{0},f) = (\bphi_{h,i}^{m})_{m=1}^{M_{i}} \in P_{0}(\mathcal{T}_{i}, \Sigma_{h,i})$ is the solution to the fully discrete subdomain problem obtained by applying the backward Euler method to \eqref{eq:disM1} on the time grid $\iT_{i}$, $i=1,2$.

\subsection{For Method 2 (GTO-Schwarz):} The two interface unknowns represent the Robin terms on each subdomain, thus we let $\zeta_{h,i} =  (\zeta_{h,i}^{m})_{m=1}^{M_{i}} \in P_{0}(\mathcal{T}_{i}, \Lambda_{h})$. The fully discrete counterpart of \eqref{eq:IP_M2} is given by
\begin{align}
&\int_{J_{1,m}} \int_{\Gamma}  \left (\zeta_{h,1} - \Pi_{12}\left (-\bphi_{h,2}(\zeta_{h,2},0,0) \cdot \bn_{1} + \alpha_{1,2}  \theta_{h,2}(\zeta_{h,2},0,0)\right )\right ) \eta_{h} \, d\gamma \, dt = \vspace{0.2cm} \nonumber\\
&\int_{J_{1,m}} \int_{\Gamma}  \Pi_{12}\left (-\bphi_{h,2}(0,f,c_{0}) \cdot \bn_{1} + \alpha_{1,2}  \theta_{h,2}(0,f,c_{0})\right ) \eta_{h} \, d\gamma \, dt, \;  \forall \eta_{h} \in \Lambda_{h}, \forall m=1, \hdots, M_{1},  \label{eq:fullIP_M2_1} \\
& \int_{J_{2,n}} \int_{\Gamma}  \left (\zeta_{h,2}-\Pi_{21}\left (-\bphi_{h,1}(\zeta_{h,1},0,0) \cdot \bn_{2} + \alpha_{2,1}  \theta_{h,1}(\zeta_{h,1},0,0)\right ) \right ) \eta_{h} \, d\gamma \, dt =\vspace{0.2cm} \nonumber\\
&\int_{J_{2,n}}  \int_{\Gamma} \Pi_{21}\left (-\bphi_{h,1}(0,f,c_{0}) \cdot \bn_{2} + \alpha_{2,1}  \theta_{h,1}(0,f,c_{0})\right ) \eta_{h} \, d\gamma \, dt, \;  \forall \eta_{h} \in \Lambda_{h},  \forall n=1, \hdots, M_{2}.  \label{eq:fullIP_M2_2}
\end{align}
where $\theta_{h,i}(\zeta_{h,i},c_{0},f) = (\theta_{h,i}^{m})_{m=1}^{M_{i}} \in P_{0}(\mathcal{T}_{i}, \Theta_{h,i})$ and $\bphi_{h,i}(\zeta_{h,i},c_{0},f) = (\bphi_{h,i}^{m})_{m=1}^{M_{i}} \in P_{0}(\mathcal{T}_{i}, \Sigma_{h,i})$ are the solution to the fully discrete subdomain problem obtained by applying the backward Euler method to \eqref{eq:disM2} on the associated time grid $\iT_{i}$, $i=1,2$.  We consider the fully discrete OSWR algorithm associated with~\eqref{eq:fullIP_M2_1}-\eqref{eq:fullIP_M2_2} using Jacobi iterations and prove that this algorithm converges. 

\subsubsection{Fully discrete OSWR algorithm with nonconforming time grids}
The OSWR algorithm reads as follows: at the $k$th iteration,  we solve, for $m=1, \hdots, M_{i}$,  the subdomain problem
\begin{equation} \label{eq:fullOSWR}
\begin{array}{l}
\left (\omega_{i} \frac{c_{h,i}^{k,m}-c_{h,i}^{k,m-1}}{\Delta t_{i,m}}, \mu_{h}\right )_{\Omega_{i}}  + \left (\nabla \cdot \bphi_{h,i}^{k,m}, \mu_{h}\right )_{\Omega_{i}}  = (f^{m},\mu_{h})_{\Omega_{i}} , \; \forall \mu_{h} \in M_{h,i}, \vspace{3pt}\\
\left (\bD_{i}^{-1} \bphi_{h,i}^{k,m}, \bv_{h}\right )_{\Omega_{i}} -\sum_{K \in \iT_{h,i}} \sum_{E \subset \partial K} u_{i,KE} \; \iU_{KE}\left (c_{i,K}^{k,m},\theta_{i,E}^{k,m}\right ) \left (\bD_{i}^{-1} \bw_{KE}, \bv_{h}\right )_{\Omega_{i}}  \\
\hspace{2cm}- \left (\nabla \cdot \bv_{h},c_{h,i}^{k,m}\right )_{\Omega_{i}}  + \sum_{K \in \iK_{h} } \left \langle \theta_{h,i}^{k,m}, \bv_{h} \cdot \bn_{K} \right \rangle_{\partial K} =0, \quad \forall \bv_{h} \in \Sigma_{h,i}, \vspace{3pt}\\
\sum_{K \in \iK_{h}} \left ( \vartheta_{h}, \bphi_{h,i}^{k,m} \cdot \bn_{K}\right )_{\partial K}  = 0, \; \forall \vartheta_{h} \in \Theta_{h}^{\Gamma, 0}, \vspace{3pt}\\
\Delta t_{i,m} \int_{\Gamma} \left (-\bphi_{h,i}^{k,m} \cdot \bn_{i}+\alpha_{i,j}  \theta_{h,i}^{k,m} \right ) \, \eta_{h} \, d\gamma \, dt \\
\hspace{3.2cm}= \int_{J_{i,m}} \int_{\Gamma} \Pi_{ij}\left (-\bphi_{h,j}^{k-1} \cdot \bn_{i} + \alpha_{i,j}  \theta_{h,j}^{k-1}\right ) \, \eta_{h} \, d\gamma \, dt, \quad \forall \eta_{h} \in \Lambda_{h},
\end{array}
\end{equation}
for $i=1,2$. 
\begin{theorem} \label{thrm:convOSWR}
Assume that $\alpha_{1,2}=\alpha_{2,1} > 0$. Algorithm~\eqref{eq:fullOSWR}, initialized by $(g_{i,j})$ in $P_{0}(\iT_{i}; \Lambda_{h})$ for $i=1,2,$ $j=3-i$, defines a unique sequence of iterates
$$ \left (c_{h,i}^{k}, \bphi_{h,i}^{k}, \theta_{h,i}^{k}\right ) \in P_{0}(\iT_{i};M_{h,i}) \times P_{0}(\iT_{i};\Sigma_{h,i}) \times P_{0}(\iT_{i};\Theta_{h,i}),
$$
that converges to the solution of the problem
\begin{equation} \label{eq:fullmulti}
\begin{array}{l}
\left (\omega_{i} \frac{c_{h,i}^{m}-c_{h,i}^{m-1}}{\Delta t_{i,m}}, \mu_{h}\right )_{\Omega_{i}}  + \left (\nabla \cdot \bphi_{h,i}^{m}, \mu_{h}\right )_{\Omega_{i}}  = (f^{m},\mu_{h})_{\Omega_{i}} , \; \forall \mu_{h} \in M_{h,i}, \vspace{3pt}\\
\left (\bD_{i}^{-1} \bphi_{h,i}^{m}, \bv_{h}\right )_{\Omega_{i}} -\sum_{K \in \iT_{h,i}} \sum_{E \subset \partial K} u_{i,KE} \; \iU_{KE}\left (c_{i,K}^{m},\theta_{i,E}^{m}\right ) \left (\bD_{i}^{-1} \bw_{KE}, \bv_{h}\right )_{\Omega_{i}}  \\
\hspace{2cm}- \left (\nabla \cdot \bv_{h},c_{h,i}^{m}\right )_{\Omega_{i}}  + \sum_{K \in \iK_{h} } \left \langle \theta_{h,i}^{m}, \bv_{h} \cdot \bn_{K} \right \rangle_{\partial K} =0, \quad \forall \bv_{h} \in \Sigma_{h,i}, \vspace{3pt}\\
\sum_{K \in \iK_{h}} \left ( \vartheta_{h}, \bphi_{h,i}^{m} \cdot \bn_{K}\right )_{\partial K}  = 0, \; \forall \vartheta_{h} \in \Theta_{h}^{\Gamma, 0}, \vspace{3pt}\\
 \Delta t_{i,m} \int_{\Gamma} \left (-\bphi_{h,i}^{m} \cdot \bn_{i}+\alpha_{i,j}  \theta_{h,i}^{m} \right ) \, \eta_{h} \, d\gamma \, dt \\
\hspace{3.2cm}= \int_{J_{i,m}} \int_{\Gamma} \Pi_{ij}\left (-\bphi_{h,j} \cdot \bn_{i} + \alpha_{i,j}  \theta_{h,j}\right ) \, \eta_{h} \, d\gamma \, dt, \quad \forall \eta_{h} \in \Lambda_{h},
\end{array}
\end{equation}
for $i=1,2$.
\end{theorem}
\begin{remark}
For simplicity, we present the convergence theorem for the case of two subdomains. However, the result still holds for a decomposition into multiple nonoverlapping subdomains $\Omega_{i}$ for $i=1, 2, \hdots, I$ and under the assumption $\alpha_{i,j}=\alpha_{j,i}$ for $i=1, 2, \hdots, I$ and $j \in \iN_{i}$, where $\iN_{i}$ denotes the set of indices of the neighbors of the subdomain $\Omega_{i}$. 
\end{remark}
\begin{proof}
As the equations are linear, we take $f=0$ and $c_{0}=0$ and prove the sequence of iterates converges to zero. We first derive an estimate for the errors $\| \theta_{h,i}\|_{0,E}$ on the edges. Following the techniques in \cite{AB85}, for $K \in \iK_{h}$ and $E \subset \partial K$ there exists a unique $\btau_{E} \in \Sigma_{h,i}$ such that $\text{supp}(\btau_{E}) \subseteq K$ and
$$ \btau_{E} \cdot \bn_{E^{\prime}} = \left \{ \begin{array}{ll} \theta_{h,i}^{k,m} & \text{if} \; E = E^{\prime},\\
0 & \text{otherwise}.
\end{array} \right .
$$
By a scaling argument, we obtain
\begin{equation} \label{eq:scaling}
h_{K} \| \btau_{E} \|_{1,K} + \| \btau \|_{0,K} \leq C h_{K}^{1/2} \| \theta_{h,i}^{k,m}\|_{0,E},
\end{equation}
where C (here and in the following) denotes a generic positive constant which is independent of the mesh size and time step size. Proceeding as in \cite{Radu14}, we take $\bv_{h}=\btau_{E}$ in \eqref{eq:fullOSWR}$_{2}$, then use \eqref{eq:scaling}, the uniform ellipticity of $\bD^{-1}$ and the uniform boundedness of $\bu$, and divide both sides by $\| \theta_{h,i}^{k,m}\|_{0,E}$ to obtain, for $E \subset \partial K$:
\begin{align*} \| \theta_{h,i}^{k,m}\|_{0,E} & \leq C \left ( h_{K}^{1/2} \| \bphi_{h,i}^{k,m}\|_{0,K} + h_{K}^{-1/2} \| c_{h,i}^{k,m}\|_{0,K} + h_{K}^{1/2} \sum_{E^{\prime} \subset \partial K} \vert E ^{\prime}\vert ( \vert c_{i,K}^{k,m} \vert + \vert \theta_{i,E^{\prime}}^{k,m}\vert \right )\\
&  \leq C \left ( h_{K}^{1/2} \| \bphi_{h,i}^{k,m}\|_{0,K} + h_{K}^{-1/2} \| c_{h,i}^{k,m}\|_{0,K} + h_{K}\sum_{E^{\prime} \subset \partial K} \|  \theta_{i,E^{\prime}}^{k,m} \|_{0,E^{\prime}}\right ).
\end{align*}
Summing over all the edges of element $K$ and for $h$ sufficiently small, we deduce that
\begin{equation} \label{eq:est_theta}
\| \theta_{h,i}^{k,m}\|_{0,E} \leq C \left ( h_{K}^{1/2} \| \bphi_{h,i}^{k,m}\|_{0,K} + h_{K}^{-1/2} \| c_{h,i}^{k,m}\|_{0,K} \right ).
\end{equation}
Next,  we choose $\mu_{h}=c_{h,i}^{k,m}$, $\bv_{h}=\bphi_{h,i}^{k,m}$ and $\mu_{h}=\theta_{h,i}^{k,m}$ in the first three equations of \eqref{eq:fullOSWR}, then add the resulting equations:
\begin{equation*}
\begin{array}{l}
\left (\omega_{i} c_{h,i}^{k,m}, c_{h,i}^{k,m}\right )_{\Omega_{i}}  - \left (\omega_{i} c_{h,i}^{k,m-1}, c_{h,i}^{k,m}\right )_{\Omega_{i}}+ \Delta t_{i,m} \left (\bD_{i}^{-1} \bphi_{h,i}^{k,m}, \bphi_{h,i}^{k,m}\right )_{\Omega_{i}} \\
-\Delta t_{i,m} \sum_{K \in \iT_{h,i}} \sum_{E \subset \partial K} u_{i,KE} \; \iU_{KE}\left (c_{i,K}^{k,m},\theta_{i,E}^{k,m}\right ) \left (\bD_{i}^{-1} \bw_{KE}, \bphi_{h,i}^{k,m}\right )_{\Omega_{i}} \\
\hspace{4cm}+ \Delta t_{i,m}\int_{\Gamma} \theta_{h,i}^{k,m} \left (\bphi_{h,i}^{k,m} \cdot \bn_{i}\right ) \, d\gamma= 0. 
\end{array}
\end{equation*}
Replacing the interface term by using the following equation
\begin{equation}
\begin{array}{ll} & \left (-\bphi_{h,i}^{k,m} \cdot \bn_{i} + \alpha_{i,j}  \theta_{h,i}^{k,m}\right )^{2}-\left (-\bphi_{h,i}^{k,m} \cdot \bn_{i} - \alpha_{j,i}  \theta_{h,i}^{k,m}\right )^{2} \\
&\hspace{1cm}= 2\left (\alpha_{i,j}+\alpha_{j,i}\right ) \theta_{h,i}^{k,m} \left (-\bphi_{h,i}^{k,m} \cdot \bn_{i}\right ) + \left (\alpha_{i,j}^{2}-\alpha_{j,i}^{2}\right ) \left ( \theta_{h,i}^{k,m} \right )^{2},
\end{array}
\end{equation}
 we obtain
\begin{equation} \label{eq:ineq1}
\begin{array}{l}
\left (\omega_{i} c_{h,i}^{k,m}, c_{h,i}^{k,m}\right )_{\Omega_{i}}  - \left (\omega_{i} c_{h,i}^{k,m-1}, c_{h,i}^{k,m}\right )_{\Omega_{i}}+ \Delta t_{i,m} \left (\bD_{i}^{-1} \bphi_{h,i}^{k,m}, \bphi_{h,i}^{k,m}\right )_{\Omega_{i}} \\
-\Delta t_{i,m}\sum_{K \in \iT_{h,i}} \sum_{E \subset \partial K} u_{i,KE} \; \iU_{KE}\left (c_{i,K}^{k,m},\theta_{i,E}^{k,m}\right ) \left (\bD_{i}^{-1} \bw_{KE}, \bphi_{h,i}^{k,m}\right )_{\Omega_{i}} \\
+ \frac{\Delta t_{i,m}}{2(\alpha_{i,j}+\alpha_{j,i})}  \int_{\Gamma}  \left (-\bphi_{h,i}^{k,m} \cdot \bn_{i} - \alpha_{j,i}  \theta_{h,i}^{k,m}\right )^{2} \, d\gamma \\
= \frac{\Delta t_{i,m}}{2(\alpha_{i,j}+\alpha_{j,i})} \int_{\Gamma}  \left (-\bphi_{h,i}^{k,m} \cdot \bn_{i} + \alpha_{i,j}  \theta_{h,i}^{k,m}\right )^{2}\, d\gamma +  \frac{\Delta t_{i,m}}{2} \int_{\Gamma}  \left ((\alpha_{j,i}-\alpha_{i,j}) \left ( \theta_{h,i}^{k,m} \right )^{2}\right )\, d\gamma.
\end{array}
\end{equation}
From the boundedness of $\bu$, \eqref{eq:est_theta} and $\vert a \vert \, \vert b \vert \leq \varepsilon a^{2} +\frac{1}{4 \varepsilon} b^{2} $ (for any $\varepsilon >0$), we have that
\begin{equation*}
\begin{array}{l}
\Delta t_{i,m}\sum_{K \in \iT_{h,i}} \sum_{E \subset \partial K} u_{i,KE} \; \iU_{KE}\left (c_{i,K}^{k,m},\theta_{i,E}^{k,m}\right ) \left (\bD_{i}^{-1} \bw_{KE}, \bphi_{h,i}^{k,m}\right )_{\Omega_{i}} \\
\leq C\Delta t_{i,m} \sum_{K \in \iT_{h,i}} \sum_{E \subset \partial K} \vert E \vert \left (\vert c_{i,K}^{k,m} \vert + \vert \theta_{i,E}^{k,m} \vert \right )  \| \bphi_{h,i}^{k,m}\|_{0,K} \\
\leq C\Delta t_{i,m} \sum_{K \in \iT_{h,i}}  \| c_{i,K}^{k,m}\|_{0,K} \, \| \bphi_{h,i}^{k,m}\|_{0,K} +
 C\Delta t_{i,m} \sum_{K \in \iT_{h,i}} h^{1/2}_{K} \sum_{E \subset \partial K} \| \theta_{i,E}^{k,m}\|_{0,E} \, \| \bphi_{h,i}^{k,m}\|_{0,K} \\
 \leq C\Delta t_{i,m} \sum_{K \in \iT_{h,i}}  \| c_{i,K}^{k,m}\|_{0,K} \, \| \bphi_{h,i}^{k,m}\|_{0,K}  + C\Delta t_{i,m} \sum_{K \in \iT_{h,i}} h_{K} \| \bphi_{h,i}^{k,m}\|_{0,K} ^{2}\\
 \leq C \varepsilon \| c_{h,i}^{k,m}\|_{0, \Omega_{i}}^{2} + \frac{C(\Delta t_{i,m})^{2}}{4\varepsilon} \| \bphi_{h,i}^{k,m}\|_{0, \Omega_{i}}^{2}+ C\Delta t_{i,m} \, h   \, \| \bphi_{h,i}^{k,m}\|_{0, \Omega_{i}} ^{2}.
\end{array}
\end{equation*}
Using this inequality, the assumptions about $\omega$, $\bD$ and $\bu$, the Cauchy-Schwarz inequality and $a^{2}-ab \geq \frac{1}{2}(a^{2}-b^{2})$, we deduce from \eqref{eq:ineq1} that
\begin{equation*}
\begin{array}{l}
\omega_{-} \left (\|c_{h,i}^{k,m}\|^{2}_{0,\Omega_{i}}-\|c_{h,i}^{k,m-1}\|^{2}_{0,\Omega_{i}}\right)+ 2\delta_{-} \Delta t_{i,m} \| \bphi_{h,i}^{k,m} \|_{0,\Omega_{i}}^{2} \\ + \frac{\Delta t_{i,m}}{(\alpha_{i,j}+\alpha_{j,i})}  \int_{\Gamma}  \left (-\bphi_{h,i}^{k,m} \cdot \bn_{i} - \alpha_{j,i}  \theta_{h,i}^{k,m}\right )^{2} \, d\gamma \\
\leq \frac{\Delta t_{i,m}}{(\alpha_{i,j}+\alpha_{j,i})}  \int_{\Gamma}  \left (-\bphi_{h,i}^{k,m} \cdot \bn_{i} + \alpha_{i,j}  \theta_{h,i}^{k,m}\right )^{2} \, d\gamma + \frac{\Delta t_{i,m}}{2} \int_{\Gamma}  \left ((\alpha_{j,i}-\alpha_{i,j}) \left ( \theta_{h,i}^{k,m} \right )^{2}\right ) d\gamma \\
+ C \varepsilon \| c_{h,i}^{k,m}\|_{0, \Omega_{i}}^{2} + \frac{C(\Delta t_{i,m})^{2}}{4\varepsilon} \| \bphi_{h,i}^{k,m}\|_{0, \Omega_{i}}^{2}+ C\Delta t_{i,m} \, h   \, \| \bphi_{h,i}^{k,m}\|_{0, \Omega_{i}} ^{2}.
\end{array}
\end{equation*}
As $c_{h,i}^{k}, \bphi_{h,i}^{k}$ and $\theta_{h,i}^{k}$ are piecewise constant on each time interval $J_{i,m}$ and as $\alpha_{i,j} =\alpha_{j,i}$ we have
\begin{equation*}
\begin{array}{l}
C_{\omega} \left (\|c_{h,i}^{k,m}\|^{2}_{0,\Omega_{i}}-\|c_{h,i}^{k,m-1}\|^{2}_{0,\Omega_{i}}\right)+ C_{\delta}\int_{ J_{i,m}} \| \bphi_{h,i}^{k} \|_{0,\Omega_{i}}^{2} \, dt \\
+ \frac{1}{(\alpha_{i,j}+\alpha_{j,i})}  \int_{ J_{i,m}} \int_{\Gamma}  \left (-\bphi_{h,i}^{k} \cdot \bn_{i} - \alpha_{j,i}  \theta_{h,i}^{k}\right )^{2} \, d\gamma \, dt \\
\leq \frac{1}{(\alpha_{i,j}+\alpha_{j,i})} \int_{ J_{i,m}} \int_{\Gamma}  \left (-\bphi_{h,i}^{k} \cdot \bn_{i} + \alpha_{i,j}  \theta_{h,i}^{k}\right )^{2} \, d\gamma \, dt,
\end{array}
\end{equation*}
where $C_{\omega} =  \omega_{-}-C \varepsilon $ and $C_{\delta}=2\delta_{-} -\frac{C\Delta t_{i,m}}{4\varepsilon} -Ch $, which are positive for sufficiently small $\varepsilon, \Delta t_{i,m}$ and $h$.
We sum over all the subintervals $J_{i,m}$ and using \eqref{eq:fullOSWR}$_4$ to obtain
\begin{equation} \label{eq:ineqT}
\begin{array}{l}
C_{\omega} \|c_{h,i}^{k,M_{i}}\|^{2}_{0,\Omega_{i}}+ C_{\delta} \int_{0}^{T} \| \bphi_{h,i}^{k} \|_{0,\Omega_{i}}^{2} \, dt + \frac{1}{(\alpha_{i,j}+\alpha_{j,i})}  \int_{0}^{T} \int_{\Gamma}  \left (-\bphi_{h,i}^{k} \cdot \bn_{i} - \alpha_{j,i}  \theta_{h,i}^{k}\right )^{2} \, d\gamma \, dt \\
\leq \frac{1}{(\alpha_{i,j}+\alpha_{j,i})} \int_{0}^{T} \int_{\Gamma}  \left (\Pi_{ij}(-\bphi_{h,j}^{k-1} \cdot \bn_{i} + \alpha_{i,j}  \theta_{h,j}^{k})\right )^{2} \, d\gamma \, dt\\
\leq \frac{1}{(\alpha_{i,j}+\alpha_{j,i})} \int_{0}^{T} \int_{\Gamma}  \left (-\bphi_{h,j}^{k-1} \cdot \bn_{i} + \alpha_{i,j}  \theta_{h,j}^{k-1}\right )^{2} \, d\gamma \, dt,
\end{array}
\end{equation}
where the last inequality is obtained due to the fact that $\Pi_{ij}$ is an $L^{2}$ projection. Define for $k \geq 1$
$$ B^{k}=\sum_{i=1}^{2}  \int_{0}^{T} \int_{\Gamma}  \left (-\bphi_{h,j}^{k} \cdot \bn_{i} + \alpha_{i,j}  \theta_{h,j}^{k}\right )^{2} \, d\gamma \, dt, \quad j=3-i,
$$
and sum \eqref{eq:ineqT} over the subdomains $i=1,2$, we deduce that
\begin{equation} \label{eq:ineqsum}
\begin{array}{l}
C_{\omega} \|c_{h,i}^{k,M_{i}}\|^{2}_{0,\Omega_{i}}+ C_{\delta} \int_{0}^{T} \| \bphi_{h,i}^{k} \|^{2}_{0,\Omega_{i}} \, dt + B^{k} \leq B^{k-1}.
\end{array}
\end{equation}
We sum over the iterates $k$ to obtain $ \|c_{h,i}^{k,M_{i}}\|^{2}_{0,\Omega_{i}}$ and $\int_{0}^{T} \| \bphi_{h,i}^{k} \|_{0,\Omega_{i}}^{2} \, dt $ converge to zero as $k \rightarrow \infty$. As  $\bphi_{h,i}^{k} \in P_{0}(\iT_{i};\Sigma_{h,i})$, the latter implies that $\| \bphi_{h,i}^{k,m} \|_{0,\Omega_{i}}^{2}$ converges to zero as $k \rightarrow \infty$ for $m=1, \hdots, M_{i}$. Finally,  it can be shown that $ \|c_{h,i}^{k,m}\|^{2}_{0,\Omega_{i}}$ converges to zero for all $m=1, \hdots, M_{i},$ and $i=1,2$ by taking $\mu_{h}=c_{h}^{m-1}$ in \eqref{eq:fullOSWR}$_{1}$.

To prove the well-posedness of a solution to \eqref{eq:fullOSWR}, it suffices to show uniqueness which is obtained by noting that \eqref{eq:ineqsum} still holds without the superscript $k$. 
\end{proof}

\begin{remark}
For the fully discrete problems with nonconforming time grids and global-in-time projections, we need to assume $\alpha_{i,j}=\alpha_{j,i}$ to perform theoretical convergence analysis. This condition is not necessary for the semi-discrete case since we can use Gronwall's lemma as in \cite{H13}. For numerical experiments, different Robin parameters still lead to convergence of the iterative algorithms. 
\end{remark}

%
%
%
%
\section{Numerical experiments}
\label{sec:NumRe}
We present three test cases to study and compare the performance of the GTP-Schur and GTO-Schwarz methods proposed in the previous sections.  In all numerical experiments, we consider $ \bD = d \pmb{I} $ isotropic and constant on each subdomain, where $ \pmb{I} $
is the 2D identity matrix.  Consequently, we may denote by $ d_{i} $, the diffusion coefficient in the subdomains.  For the first two test cases, we consider a decomposition into two subdomains. In Test case~1,  the same constant coefficients are used in the subdomains, while in Test case~2, discontinuous coefficients are used with different Pecl\'et numbers to check the robustness of the methods when advection is dominant.  Test case~3 is a prototype for the simulation of the transport around a surface nuclear waste storage where the geometry of the computational domain is complex and the physical coefficients are highly variable.  The domain is decomposed into six subdomains and time windows are used for long time simulations.  

We aim to investigate the accuracy and the convergence speed of the iterative algorithms. Particularly,  we will verify the performance of Neumann-Neumann preconditioner for GTP-Schur and optimized parameters for GTO-Schwarz.  For GTP-Schur,  we use the following
formula for calculating the weights $\sigma_{i,j}$ in~\eqref{eq:precondS} (see~\cite{MB96,H13}):
\begin{equation*}
\sigma_{i,j} = \left (\frac{d_{i}}{d_{i}+d_{j}}\right )^{2}, \; i =1, 2, \hdots, I, \; j \in \iN_{j}.
\end{equation*}
For GTO-Schwarz,  we use two-sided optimized Robin parameters, i.e. $\alpha_{i,j} \neq \alpha_{j,i}$, obtained by numerically minizing the continuous convergence factor of the OSWR algorithm \cite{GHK07,Japhet12}.

%
%
%
%
	\subsection{Test case 1: with a known analytical solution}
	
%
%
We first verify the accuracy in space and in time of the proposed algorithms by considering a test case with the exact solution is given by $$ u(x,y,t)=e^{-4t} \, \sin (\pi x) \, \sin (\pi y),$$ on the unit square $\Omega=(0,1)^{2}$. 
We split $ \Omega $ into two nonoverlapping subdomains $ \Omega_{1} = (0,0.5) \times (0,1) $ and $ \Omega_{2} = (0.5,1) \times (0,1) $. Constant parameters are imposed on the whole domain: $ \omega_{i}=1 $, $ \bu_{i} =(1, \; 1)^{T} $, $d_{i}=1$, for $i=1,2$. For the spatial discretization, we consider a conforming rectangular mesh with size $ h_{1}=h_{2}= h $. For the time discretization,  we use nonconforming time grids with $ \Delta t_{1} \neq \Delta t_{2}$. The interface problem associated with each method is solved iteratively by GMRES with a zero initial guess on the interface; the iteration stops when the relative residual is smaller than $10^{-6}$. 

%
%
In Table~\ref{tab:test1accuSpace}, we show the relative $L^{2}(\Omega)-$norm errors of $c$ and $\bphi$ at $T=0.1$ with fixed time step sizes $\Delta t_{1}=T/80$ and $\Delta t_{2}=T/60$ and a decreasing mesh size $h$.  The number of subdomain solves are also reported for the global-in-time Schur (GT-Schur) method with or without Neumann-Neumann preconditioner and the GTO-Schwarz method with optimized parameters.  Note that one iteration of the GT-Schur method with the preconditioner costs twice as much as one iteration of
the GTO-Schwarz method (in terms of number of subdomain solves),  thus we show the number of subdomain solves (instead of number of iterations) to compare the two methods.  We observe that the errors by both methods are almost identical, and the order of accuracy in space is preserved with nonconforming time grids. For GT-Schur, the preconditioner significantly improves the convergence speed,  and the convergence is almost independent of $h$.  GTO-Schwarz converges a little slower than GT-Schur for this case,  and the number of GMRES iterations increases slightly when $h$ is decreasing. Note that for steady-state problems \cite{Gander06},  the convergence factor of the optimized Schwarz algorithm (i.e. Jacobi iterations) behaves like $1-O(h^{1/4})$. 
%
%
\begin{table}[h] 
\setlength{\extrarowheight}{2pt}
 \scriptsize
	\begin{tabular}{| l | c | c | c | c | c | c | c |} \hline 
		\multirow{3}{*}{$h$} & \multicolumn{4}{c|}{Method 1: GT-Schur} & \multicolumn{3}{c|}{Method 2: GTO-Schwarz}  \\ \cline{2-8}
		& \multicolumn{2}{c|}{$L^{2}$ errors} & \multicolumn{2}{c|}{\# subdomain solves} & \multicolumn{2}{c|}{$L^{2}$ errors} & \# subdomain \\ \cline{2-3} \cline{4-7}
		& $c$ & $\bphi$ & {\tiny No Precond.} & {\tiny With Precond.} & $c$ & $\bphi$ & {solves}\\ \hline
		$1/20$ &  0.0641	\phantom{[1.00]}& 0.0453	\phantom{[1.00]}& 29  & 	12 & 0.0641 \phantom{[1.00]}& 0.0454 \phantom{[1.00]}& 16 \\
		$1/40$ & 0.0321	[1.00] & 0.0227	[1.00]& 39  & 12 & 0.0321 [1.00]	& 0.0227 [1.00] & 16\\
		$1/80$ &  0.0160	[1.00] & 0.0114  [0.99]  & 54 & 12 &  0.0160 [1.00]  & 0.0114 [0.99] & 20 \\
		$1/160$ & 0.0080	[1.00] & 0.0057 [1.00] & 76 & 14  & 0.0080 [1.00] & 0.0058 [0.98] & 22  \\ \hline
	\end{tabular}  
	\caption{[Test case 1] Accuracy in space, the convergence rates are shown in square brackets.} 
\label{tab:test1accuSpace}
\end{table}
%
%
%
%

In Table~\ref{tab:test1accuTime}, we show the relative $L^{2}(\Omega)-$norm errors of $c$ and $\bphi$ at $T=1$ with fixed $h=1/200$ and decreasing time step sizes $\Delta t_{1} =\sfrac{3}{4}\Delta t_{2}$.  The accuracy in time is preserved with nonconforming time grids.  In addition,  the errors obtained by the two methods are not the same, especially when the time step sizes are large.  This is due to the use of different projection operators, which makes the two methods yield different solutions at convergence (cf. Section~\ref{sec:time}).  However, as $\Delta t_{i}, \; i=1,2,$ become smaller,  numerical results suggest that the two methods converge to the same continuous-in-time solution.  Again, the preconditioner as well as optimized parameters work well in terms of convergence speed. 
%
%
\begin{table}[h] 
\setlength{\extrarowheight}{2pt}
\scriptsize
	\begin{tabular}{| l | c | c | c | c | c | c | c |} \hline 
		 \multirow{3}{*}{$\Delta t_{2}$} & \multicolumn{4}{c|}{Method 1: GT- Schur} & \multicolumn{3}{c|}{Method 2: GTO-Schwarz}  \\ \cline{2-8}
		& \multicolumn{2}{c|}{$L^{2}$ errors} & \multicolumn{2}{c|}{\# subdomain solves} & \multicolumn{2}{c|}{$L^{2}$ errors} & \# subdomain \\ \cline{2-3} \cline{4-7}
		& $c$ & $\bphi$ & {\tiny No Precond.} & {\tiny With Precond.}  & $c$ & $\bphi$ & {solves}\\ \hline
		$T/6$ &  0.1186	\phantom{[1.00]}& 0.1315	\phantom{[1.00]}& 83 & 14  & 0.2524 \phantom{[1.00]}& 0.2712 \phantom{[1.00]}& 18\\ 
		$T/12$ & 0.0520	 [1.19]	& 0.0579	[1.18]	& 83 & 14 & 0.0922  [1.45]& 0.1042 [1.38]& 18\\
		$T/24$ &  0.0251 	[1.05]	& 0.0277 [1.06]	& 93 & 14 & 0.0369 [1.32]& 0.0422 [1.30]& 20 \\
		$T/48$ & 0.0134 	[0.91]	& 0.0141 [0.97]		& 104 & 14 & 0.0160 [1.21]& 0.0173 [1.29]& 24 \\ \hline 
	\end{tabular}  
	\caption{[Test case 1] Accuracy in time, the convergence rates are shown in square brackets.} 
\label{tab:test1accuTime}
\end{table}
%
%

\subsection{Test case 2: with piecewise discontinuous coefficients}
Next, we analyze the convergence of the iterative algorithms. Towards that end, we consider the error equation with the same two nonoverlapping subdomains as in Test case~1. The porosity is $ \omega _{1} = \omega_{2} =\omega = 1 $. The diffusion and advection coefficients, $ d_{i} $ and $ \bu_{i} $ for $ i~=1,2, $ are constant in each subdomain and  discontinuous across the interface. Their values are given in Table~\ref{tab:Test2Discont} for the diffusion dominant, mixed regime and advection dominant problems, respectively. Note that the global P\'eclet number in each subdomain is computed by $$ \text{Pe}_{G,i}:=\frac{H\mid \bu_{i} \mid}{d_{i}}, \; i=1,2,$$ where H is the size of the domain (in this case, $ H=1 $). 
In space, we use a conforming rectangular mesh $h= 1/100 $; in time, nonconforming time grids are considered with $\Delta t_{1}~=\sfrac{3}{4} \Delta t_{2}$ and $\Delta t_{2}=1/75$.
\begin{table}[h]
\footnotesize
\centering
\begin{tabular}{|l|l|l|l|l|l|l|}
  \hline
  		Problems							&$  d_{1}$	 		& $ \bu_{1} $		& $ \text{Pe}_{G,1}$ &$  d_{2}$	 & $ \bu_{2} $				& $ \text{Pe}_{G,2}$ 	 \\ \hline
   (a) Diffusion dominance    &  $1$  &   	$(-0.02, \; -0.5)^{T}$ &  $\approx 0.5$  & $0.1$ & $(-0.02, \; -0.05)^{T}$ & $\approx 0.5$	\\ \hline 
  (b) Mixed regime     &  $0.01$  &   	$(-0.02, \; -0.5)^{T}$ &  $\approx 50$  & $0.1$ & $(-0.02, \; -0.05)^{T}$ & $\approx 0.5$		\\ \hline 
  (c) Advection dominance     &  $0.02$  &   	$(0.5, \; 1)^{T}$ &  $\approx 56$  & $0.002$ & $(0.5, \; 0.1)^{T}$ & $\approx 255$		\\ \hline 
\end{tabular}
\caption{[Test case 2] Data for the discontinuous test case.}
\label{tab:Test2Discont} \vspace{-0.2cm}
\end{table}

\begin{figure}[!ht]
\centering
    \begin{subfigure}[b]{0.34\textwidth}
    \centering
        \includegraphics[width=\textwidth]{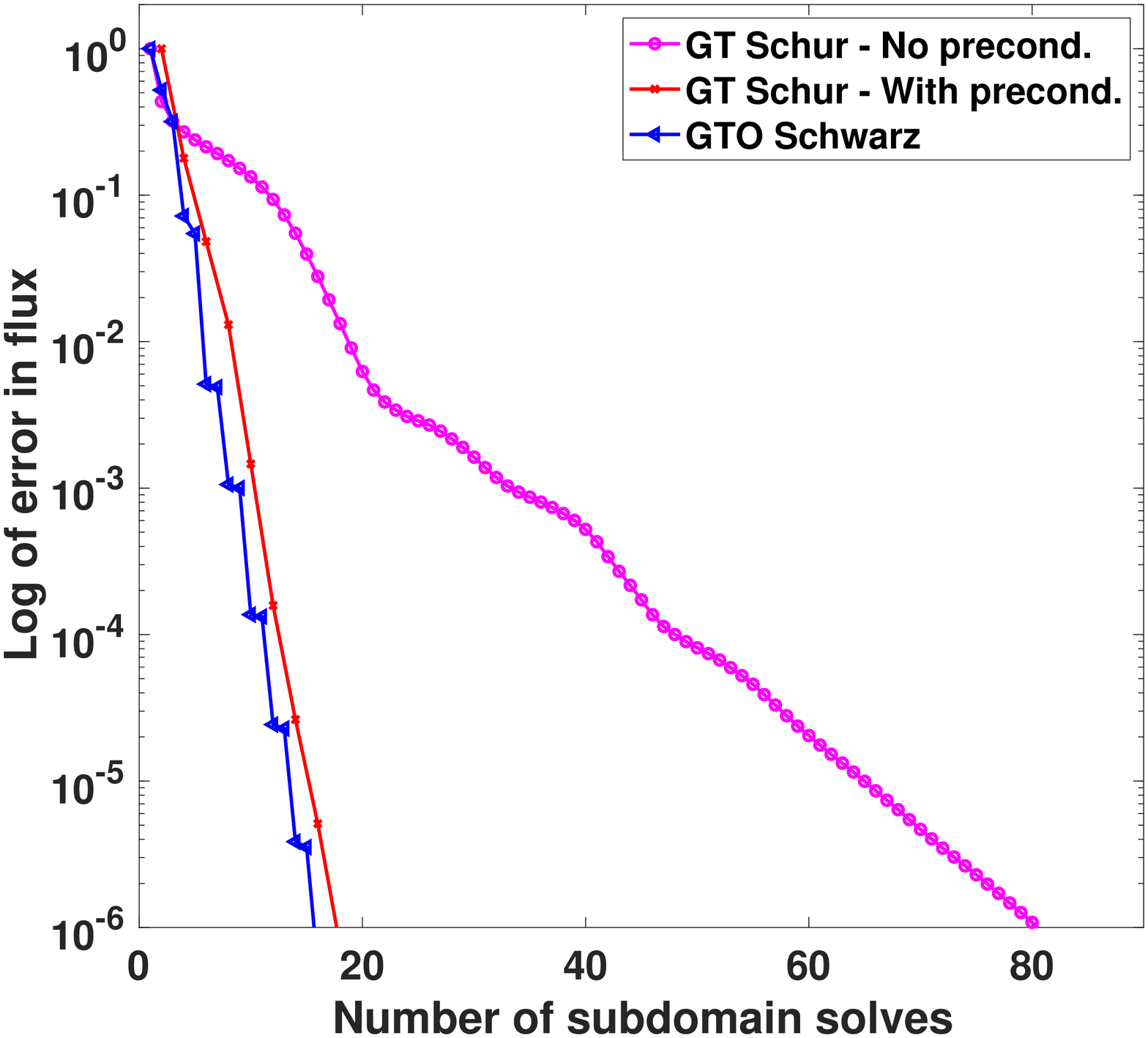}
        \caption{Diffusion dominance \\ $\max Pe_{G,i} \approx 0.5$}
    \end{subfigure}\hspace{-0.25cm}
    \begin{subfigure}[b]{0.34\textwidth}
     \centering
        \includegraphics[width=\textwidth]{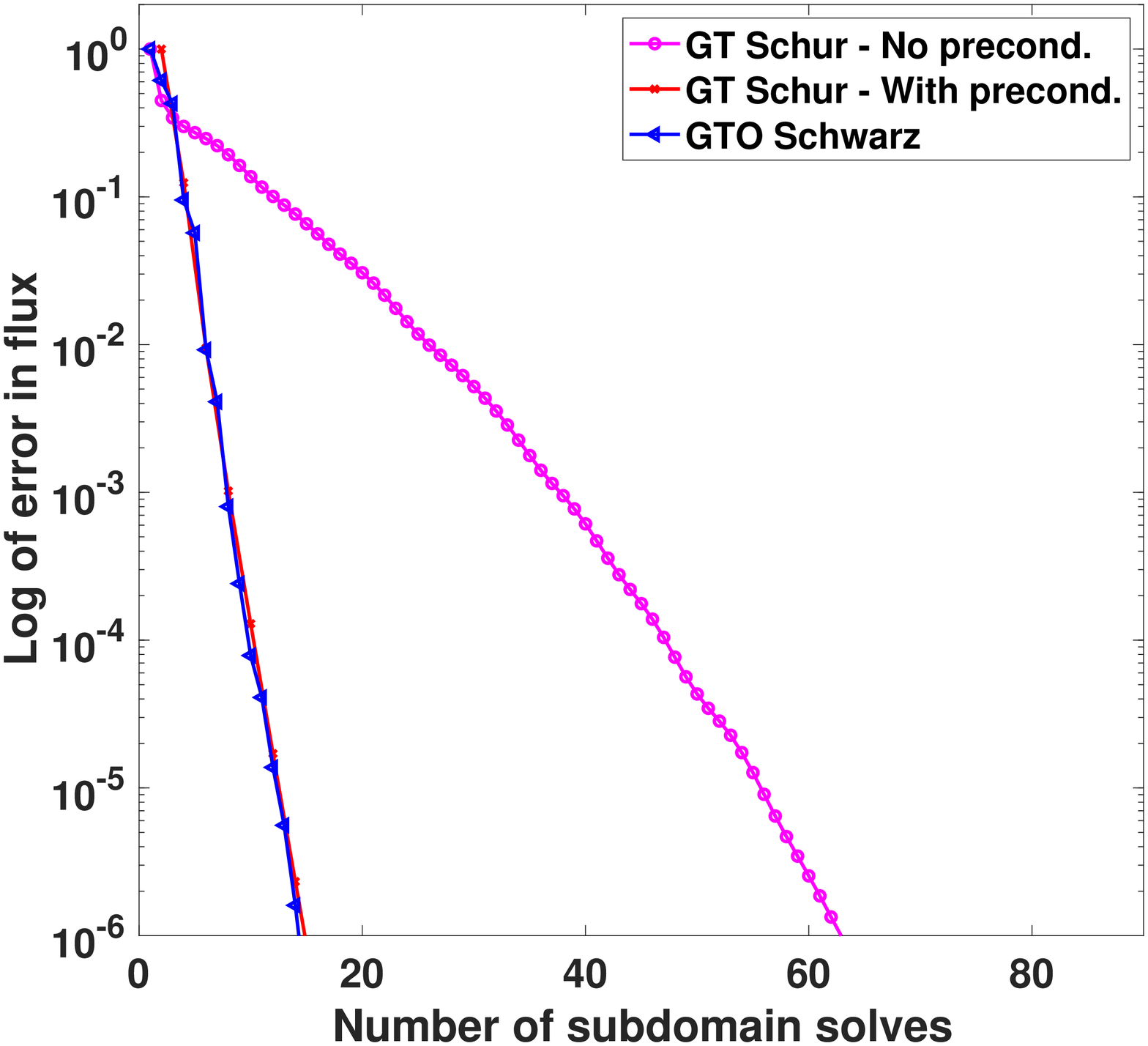}
        \caption{Mixed regime \\ $\max Pe_{G,i} \approx 50$}
    \end{subfigure}\hspace{-0.25cm}
    \begin{subfigure}[b]{0.335\textwidth}
     \centering
        \includegraphics[width=\textwidth]{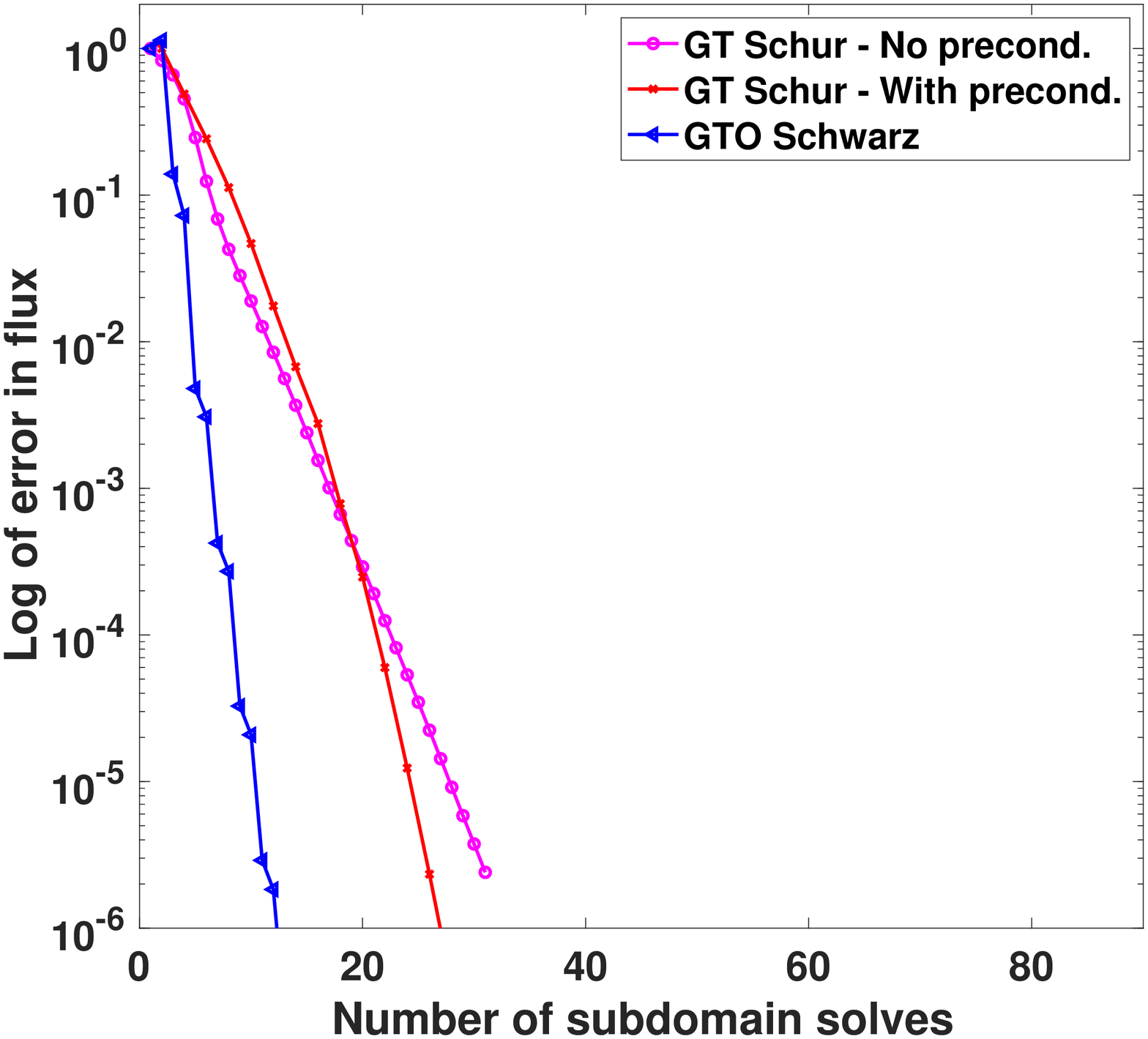}
        \caption{Advection dominance \\ $\max Pe_{G,i} \approx 255$}
    \end{subfigure} \vspace{-0.2cm}
    \caption{[Test case 2] Convergence curves by GMRES for different P\'eclet numbers: $L^{2}-$norm errors in the flux $\bphi$ at $T=1$ for Method 1 (GT-Schur) with no preconditioner (magenta curves) and with the Neumann-Neumann preconditioner (red curves), and Method 2 (GTO-Schwarz) (blue curves).}\label{fig:test2GMRESconv} 
\end{figure}

Figure~\ref{fig:test2GMRESconv} show the errors (in logarithmic scale) in $L^{2}(\Omega)-$norm of the flux variable $\bphi$ versus the number of subdomain solves using GMRES with a random initial guess (similar convergence curves are obtained for the scalar variable $c$, and are omitted).  Three algorithms are considered: GT-Schur with no preconditioner (magenta, circle),  GT-Schur
with the preconditioner (red, x-mark) and GTO-Schwarz (blue, triangle). We observe that
for GT-Schur, the preconditioner works well in the case the Pecl\'et number is not so large. (i.e.  Problems (a) and (b) in Table~\ref{tab:Test2Discont}).  When the Pecl\'et number is sufficiently large (i.e.  Problem (c)),  the convergence of GT-Schur with or without preconditioner is quite the same.  For GTO-Schwarz, the convergence
speed does not significantly change with the Pecl\'et number.  GT-Schur with the preconditioner is comparable with GTO-Schwarz when diffusion is dominant. When advection is dominant, GTO-Schwarz converges faster than GT-Schur with or without preconditioner (at least by a factor of 2.17).  We remark that when operator splitting is used as in \cite{H17, Hthesis},  the GT-Schur approach with the preconditioner converges even slower than without preconditioner when advection is dominant (cf.  Figure 3 in~\cite{H17} and Figure 3.8 in \cite[Chapter 3]{Hthesis} where the Pecl\'et numbers are approximately $100$ and $100\sqrt{2}$ respectively). 

\begin{figure}[!ht]
\centering
    \begin{subfigure}[b]{0.49\textwidth}
    \centering
        \includegraphics[width=\textwidth]{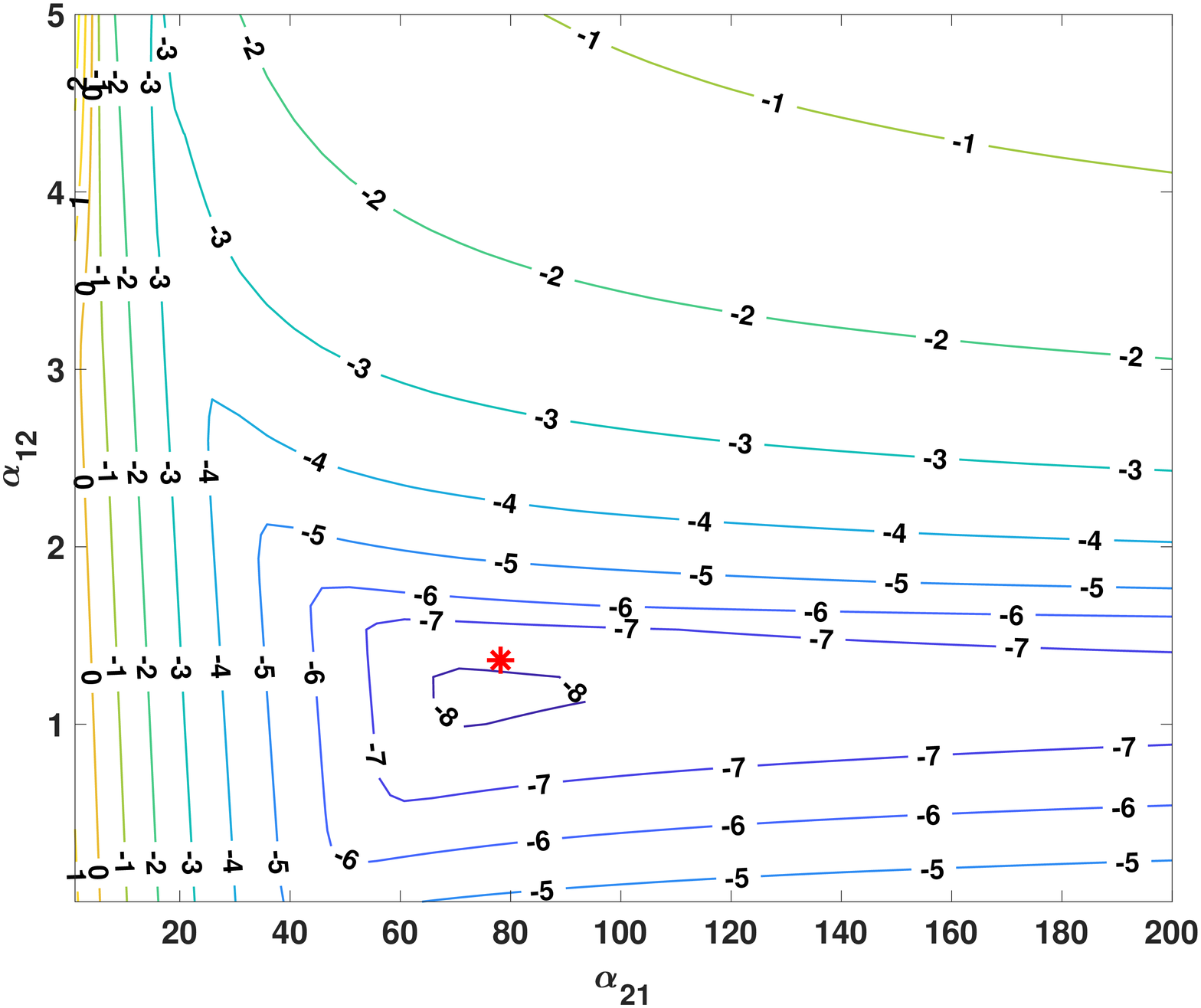}
        \caption{Diffusion dominance \\ 25 Jacobi iterations}
    \end{subfigure}\hspace{-0.25cm}
    \begin{subfigure}[b]{0.49\textwidth}
     \centering
        \includegraphics[width=\textwidth]{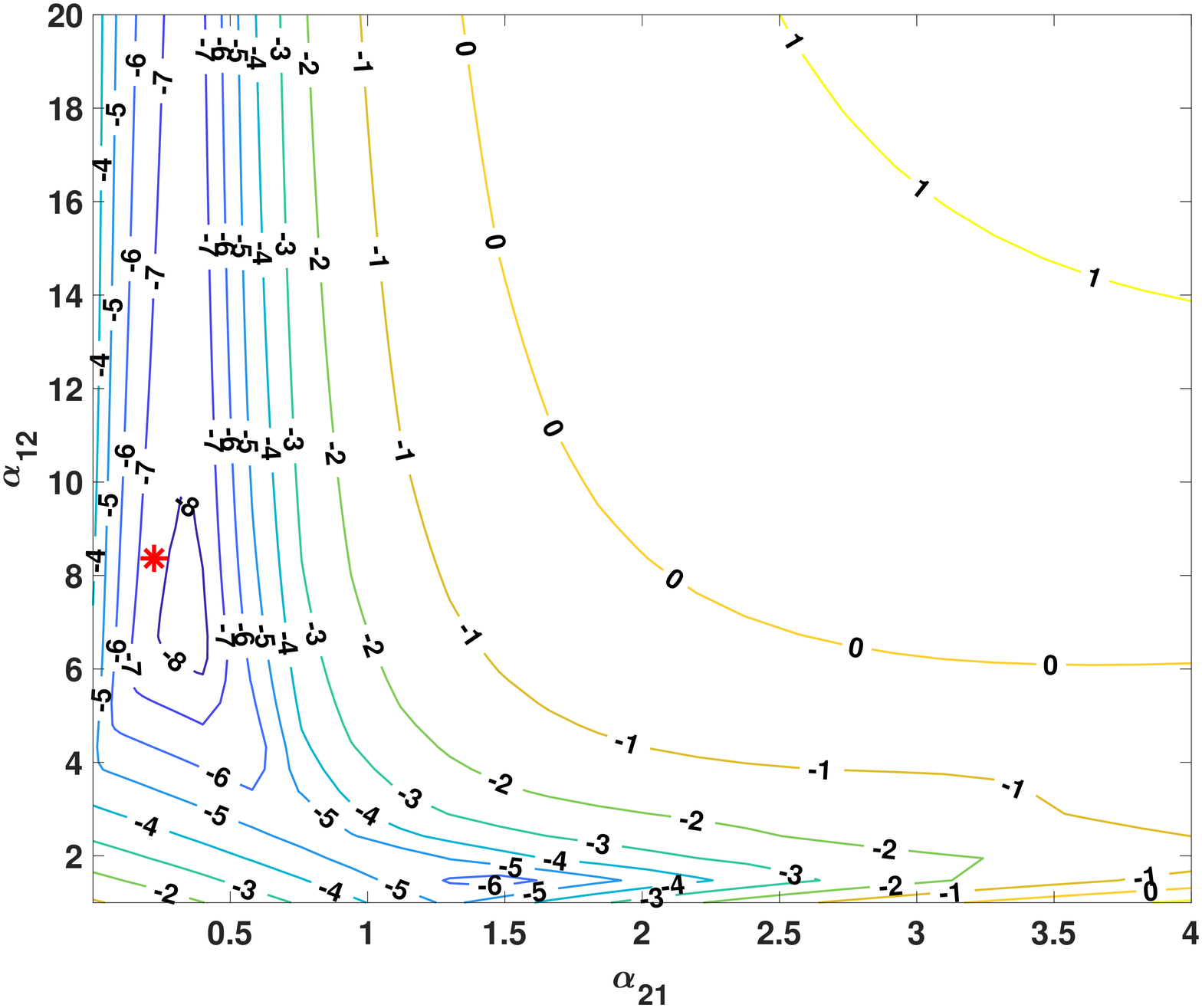}
        \caption{Mixed regime \\ 25 Jacobi iterations}
    \end{subfigure}\hspace{-0.25cm}
    \begin{subfigure}[b]{0.49\textwidth}
     \centering
        \includegraphics[width=\textwidth]{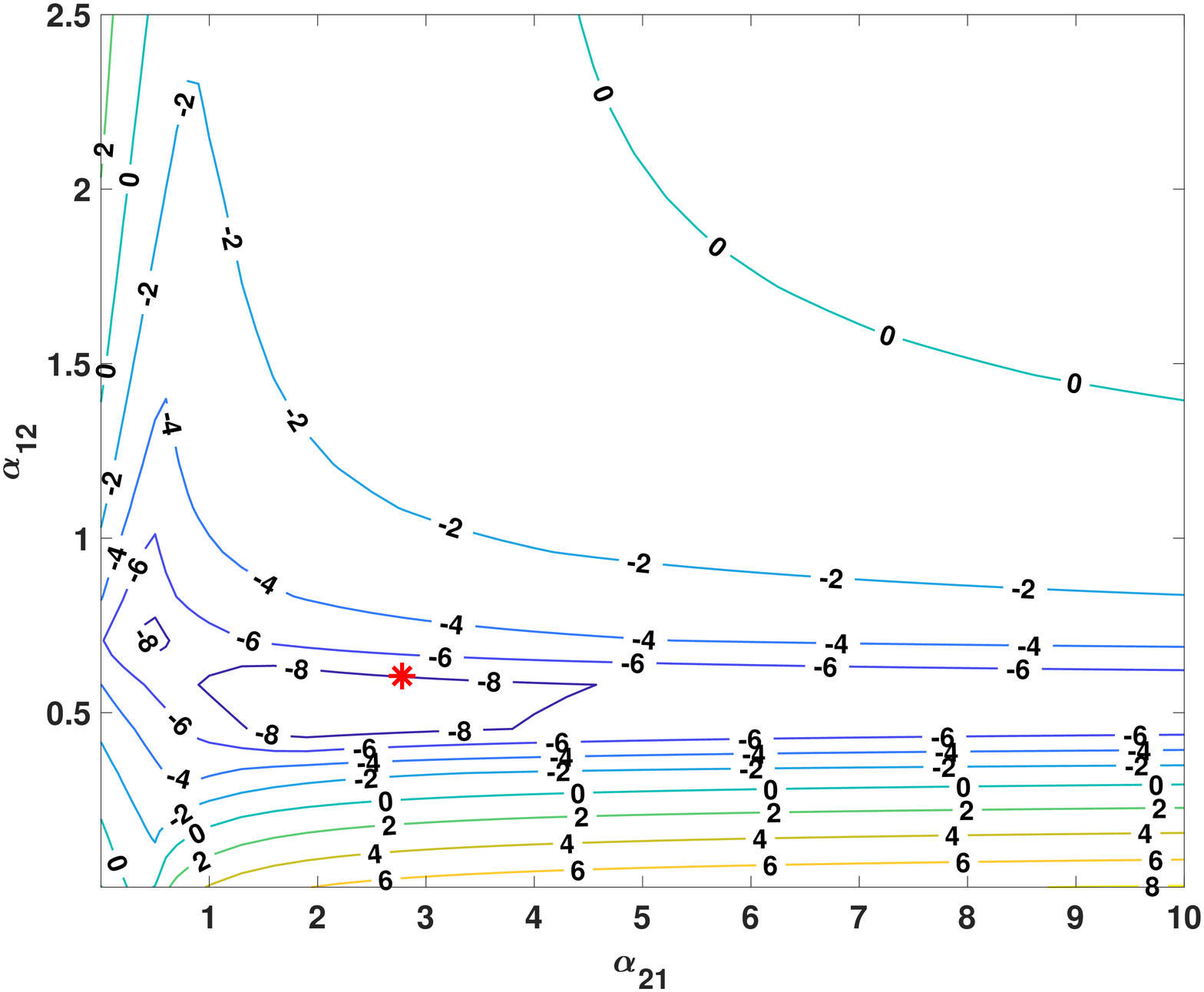}
        \caption{Advection dominance \\ 20 Jacobi iterations}
    \end{subfigure} \vspace{-0.2cm}
    \caption{[Test case 2] Level curves for the relative residuals (in logarithmic scales) after a fixed number of Jacobi iterations for various values of $\alpha_{1,2}$ and $\alpha_{2,1}$. The red star shows the optimized values computed by numerically minimizing the continuous convergence factor of the OSWR algorithm. }\label{fig:test2optpara} 
\end{figure}

To verify the performance of the optimized parameters, we show in Figure~\ref{fig:test2optpara} the relative residuals (in logarithmic scale) for various values of the parameters $\alpha_{1,2}$ and $\alpha_{2,1}$ after a fixed number of Jacobi iterations.  We see that for different sets of parameters,  the pair of optimized Robin parameters (red star) is located close to those giving the smallest relative residual after the same number of iterations.

\begin{figure}[!ht]
\centering
    \begin{subfigure}[b]{0.48\textwidth}
    \centering
        \includegraphics[width=\textwidth]{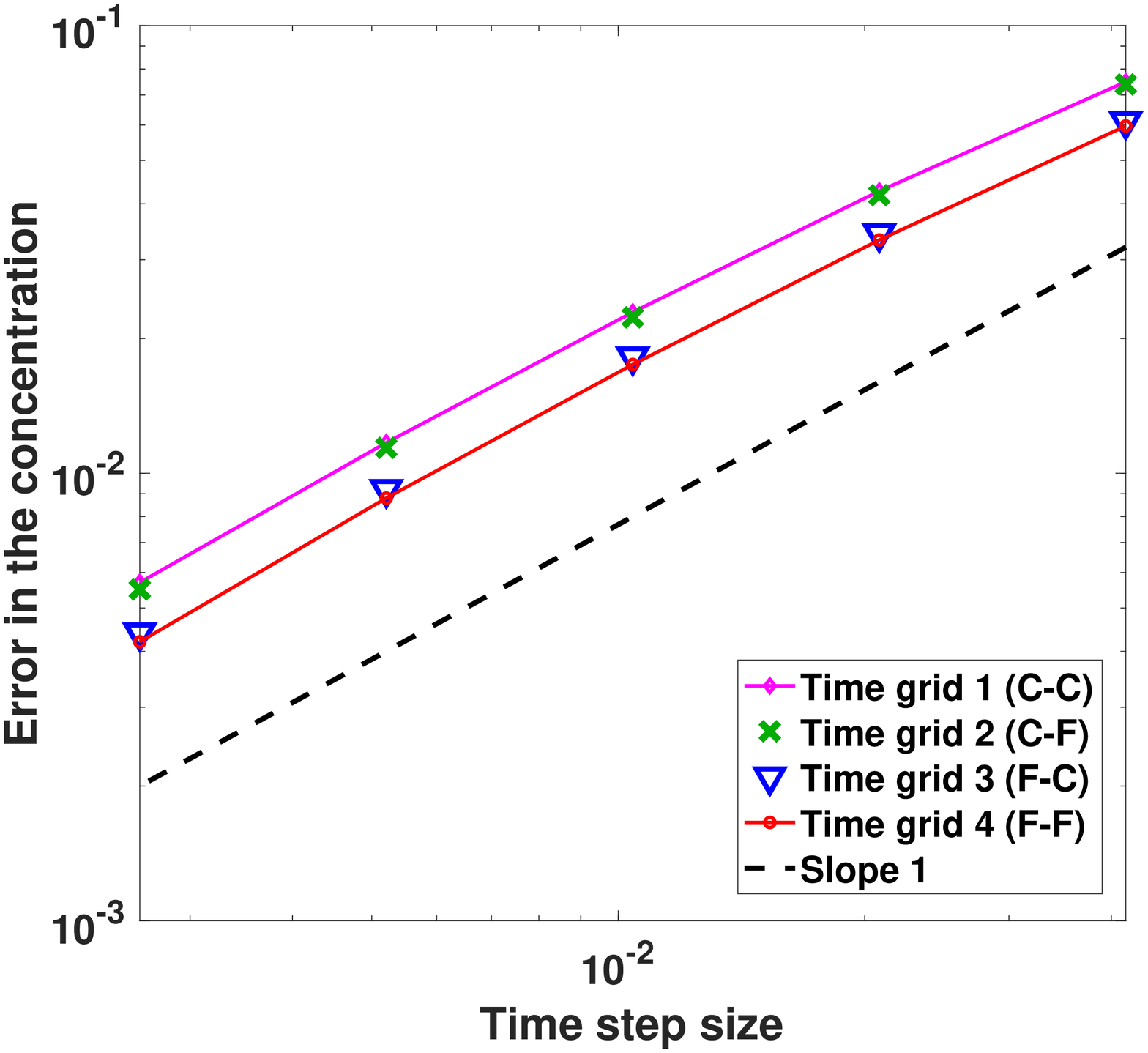}
    \end{subfigure} \hspace{4pt}
    \begin{subfigure}[b]{0.48\textwidth}
     \centering
        \includegraphics[width=\textwidth]{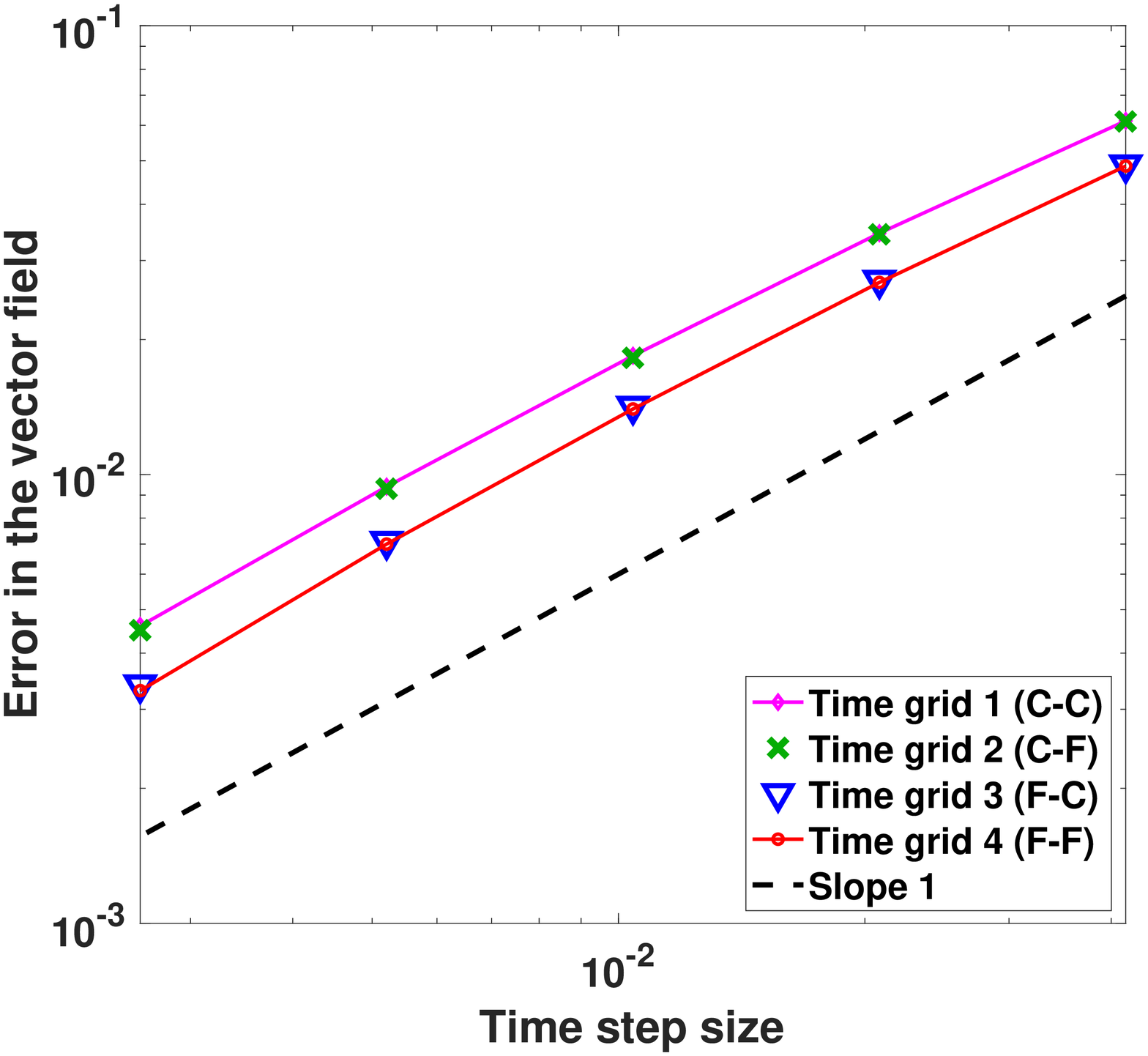}
    \end{subfigure} \vspace{-0.2cm}
    \caption{[Test case 2] Errors in the concentration $c$ (left) and the vector field $\bphi$ (right) between the reference and multidomain solutions.}\label{fig:test2accuracytime} 
\end{figure}

We now investigate whether the nonconforming time grids preserve the accuracy in time.  We consider the  advection dominant problem (i.e. Problem (c)) with the same coefficients given in Table~\ref{tab:Test2Discont}. Homogeneous Dirichlet conditions are imposed on the boundary,  the source term is $ f(x,y, t)=\exp(-100((x-0.2)^2+(y-0.2)^2)),$ and the initial condition $p_{0}(x,y)=xy(1-x)(1-y)\exp(-100((x-0.2)^2+(y-0.2)^2)).$ We use four initial time grids with $ \Delta t_{c}=T/12$ and $ \Delta t_{f}=T/16 $ where $T=0.5$:
\begin{itemize} \itemsep0pt
	\item Time grid 1 (coarse-coarse): conforming with $ \Delta t_{1} = \Delta t_{2} = \Delta t_{c} $.
	\item Time grid 2 (coarse-fine): nonconforming with $ \Delta t_{1} = \Delta t_{c} $ and
              $ \Delta t_{2} = \Delta t_{f} $.
	\item Time grid 3 (fine-coarse): nonconforming with $ \Delta t_{1} = \Delta t_{f} $ and
              $ \Delta t_{2} = \Delta t_{c} $.
	\item Time grid 4 (fine-fine): conforming with $ \Delta t_{1} = \Delta t_{2} = \Delta t_{f} $.
\end{itemize}
The time steps are then refined several times by a factor of 2. In space, we fix a conforming rectangular mesh with $h=1/200$, and we compute a reference solution by solving problem~\eqref{eq:dismono} directly
on a very fine time grid, with $ \Delta t = \Delta t_{f}/ 2^{7} $.
The converged DD solution is such that the relative residual is smaller than $ 10^{-8} $.  We show in Figure~\ref{fig:test2accuracytime} the relative errors at $T=0.5$ versus the time step $ \Delta t = \max (\Delta t_{c}, \Delta t_{f}) $. We only give the results for GTP-Schur because the curves for GTO-Schwarz look exactly the same.  We observe that first order convergence is preserved in the nonconforming case.  The errors obtained in the
nonconforming case with a fine time step in $\Omega_{1}$ where the parameters are large (Time grid 3 with blue triangle markers) are nearly the same as in the finer
conforming case (Time grid 4, in red with circle markers).  On the other hand,  the errors obtained in the
nonconforming case with a fine time step in $\Omega_{2}$ where the parameters are small (Time grid 2 with green x-markers) are close to those by the coarse
conforming case (Time grid 1, in magenta with diamond markers).  Thus using nonconforming grids can
adapt the time steps in the subdomains depending on the physical parameters and limit the computational cost locally,  while preserving almost the same accuracy as in the finer conforming case.

	\subsection{Test case 3: A simulation for a surface, nuclear waste storage}
	\label{A2subsec:ANDRA}
\begin{figure}[!http]
\centering
\begin{minipage}[c]{0.9 \linewidth}
\begin{center}
\includegraphics[scale=0.5]{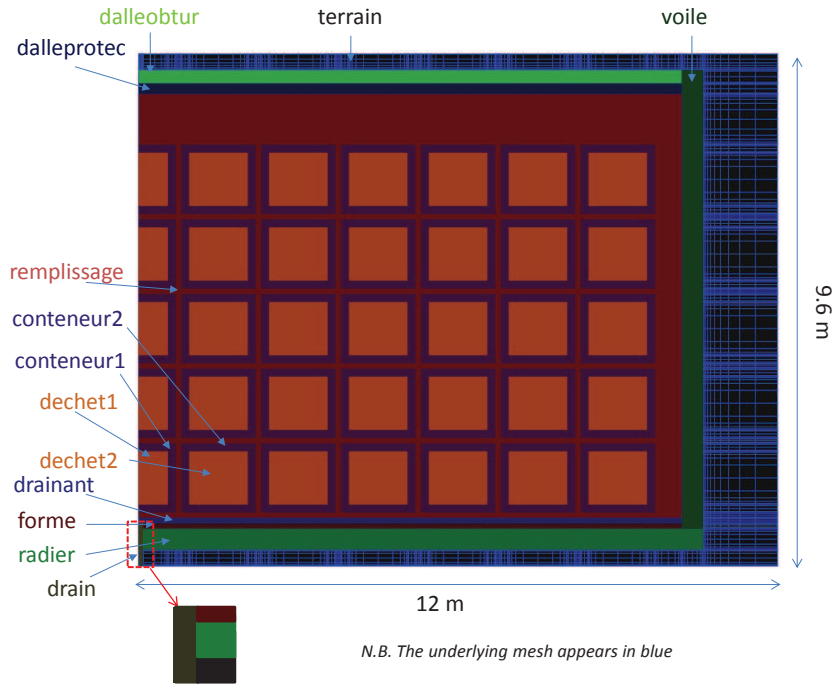} 
\end{center}
\end{minipage}  \vspace{-0.2cm}
\caption{[Test case 3] The geometry of the test case \cite{H17}.} 	
\label{A2Fig:ANDRAdomain} \vspace{-0.2cm}
\end{figure}

\begin{figure}[!http]
\vspace{0.2cm}
\begin{minipage}{0.4 \linewidth}
\hspace{-0.6cm}\includegraphics[scale=0.29]{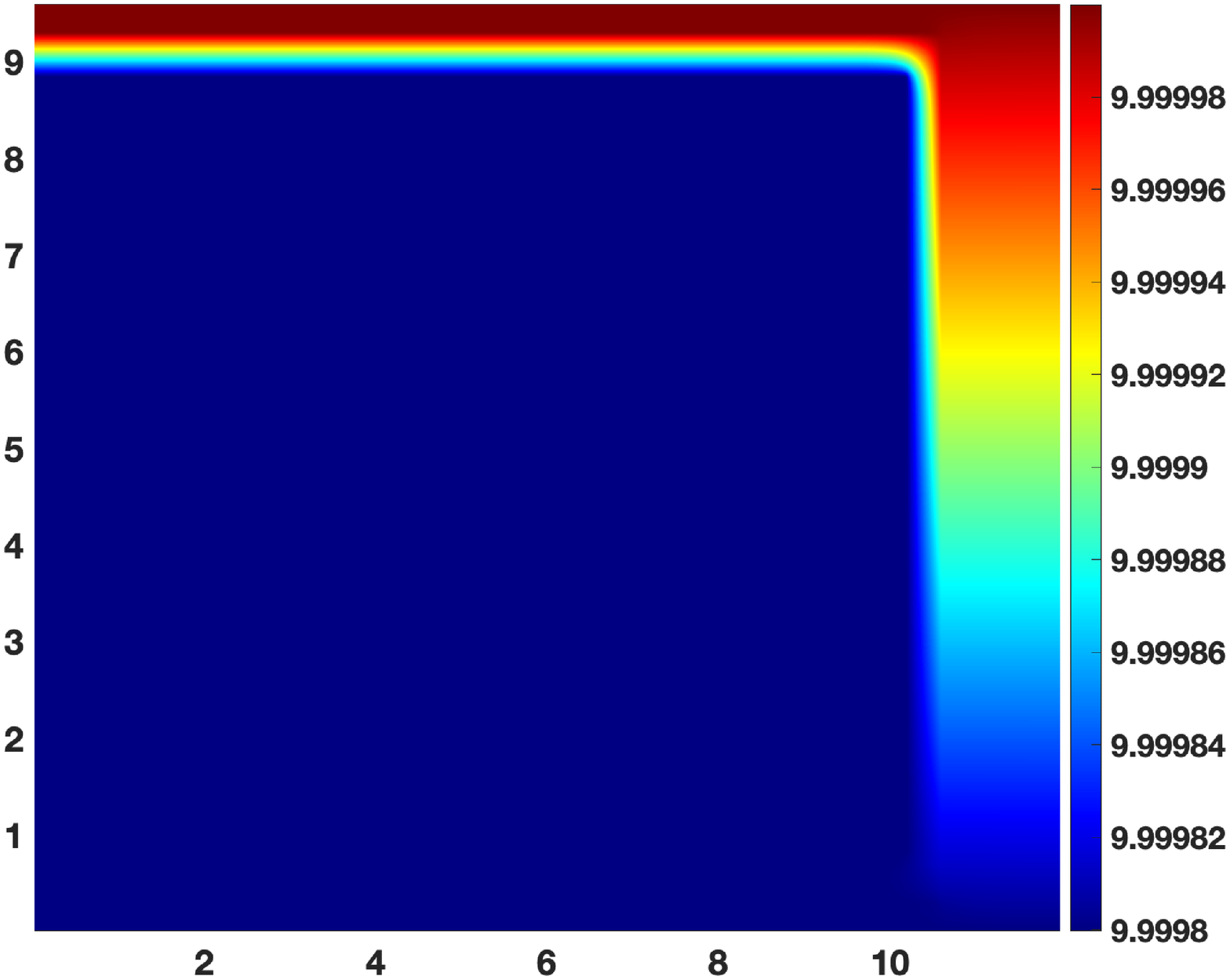} 
\end{minipage} \hspace{2pt}
\begin{minipage}{0.5 \linewidth}
\vspace{-0.1cm}
\setlength{\unitlength}{1pt} 
\begin{picture}(140,100)(0,0)
\thicklines
\put(45,-40){\includegraphics[height=6.7cm]{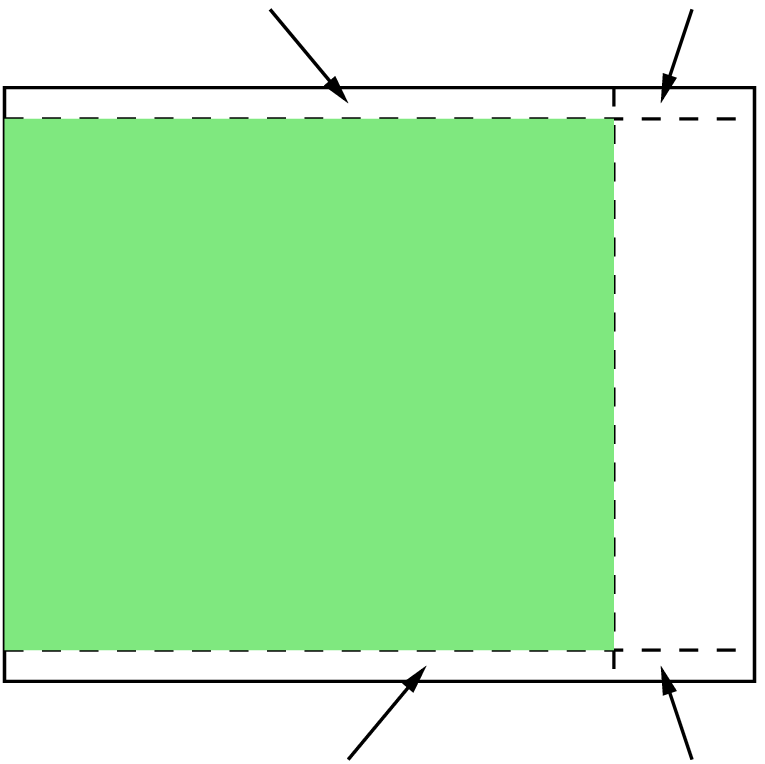} \\}
\put(120,-37){$ \Omega_{1} $}
\put(222,-37){$ \Omega_{2} $}
\put(120,45){$ \Omega_{3}$}
\put(210,45){$ \Omega_{4} $}
\put(100,145){$ \Omega_{5} $}
\put(220,145){$ \Omega_{6}$}
\end{picture}
\end{minipage} 
\caption{[Test case 3] The hydraulic head field and the decomposition of the domain.} 	
\label{A2Fig:ANDRADarcyDecom}   \vspace{-0.2cm}
\end{figure}

%
%
Finally, we consider a test case introduced in \cite{H17} and designed by ANDRA\footnote{The French agency for nuclear waste management} as a protype for simulating a surface storage of short half-life nuclear waste.  The computational domain is depicted in Figure~\ref{A2Fig:ANDRAdomain} with different physical zones, where the waste is stored in square boxes ("dechet" zone). The properties of these zones are given in Table~\ref{A2Tab:ANDRAdata}. Note that in our calculation, we use the effective diffusion, defined by $d_{\text{eff}} = \omega \times d_{\text{m}} $. The advection field is governed by Darcy's law together with the law of mass conservation:
\begin{equation} \label{A2ANDRADarcyflow}
\begin{array}{rll} \bu & = - K \nabla \mathfrak{h} & \text{in} \; \Omega, \\
\Div \bu & = 0 & \text{in} \; \Omega,
\end{array} 
\end{equation}
where $h$ is the hydraulic head field, $\bu$ is the Darcy velocity and $K$ is the hydraulic conductivity. Dirichlet conditions are imposed on top, $\mathfrak{h}= 10$m and on bottom $ \mathfrak{h} = 9.998$m of the domain and no flow boundary on the left and right sides.  For the transport problem,  the final time is $T_{f}=500$ years,  the source term is $ f = 0 $ and the initial condition is such that 
\begin{equation*} \label{A2ANDRA:IC}
c_{0}= \left \{ \begin{array}{ll} 1, & \text{in "dechet1" and "dechet2"}, \\
0, & \text{elsewhere}.
\end{array} \right .
\end{equation*}
Boundary conditions of the transport problem are homogeneous Dirichlet conditions on top and bottom, and homogeneous Neumann conditions on the left and right hand sides.  
\begin{table}[!http]
\centering
\begin{tabular}{|l|l|l|l|l|}
  \hline
	    		Zone               						& Hydraulic conductivity       & Porosity                    & Molecular diffusion  				 \\ 
	    													& $ K $ (m/year) 				  & $ \omega $					& $ d_{\text{m}}$ (m$^{2}$/year)   \\ \hline
	    		terrain								    & $ 94608 $   									& $ 0.30 $                  & $ 1 $    	 					\\ \hline 
	    		radier									& $ 3.1536 \, 10^{-4} $   					& $ 0.15 $                  & $ 6.31 \, 10^{-5}$     	  \\ \hline 	
	    		forme									& $ 3.1536 \, 10^{-3} $   					& $ 0.20 $                  & $ 1.58 \, 10^{-3} $     	  \\ \hline 
	    		drainant 								& $ 94608 $   									& $ 0.30 $                  & $ 5.36 \, 10^{-2} $     	  \\ \hline
	    		voile 									& $ 3.1536 \, 10^{-3} $   					& $ 0.20 $                  & $ 1.58 \, 10^{-3} $   	  \\ \hline
	    		remplissage 							& $ 5045.76 $   								& $ 0.30 $                  & $ 5.36 \, 10^{-2} $     	  \\ \hline
	    		dalleprotec 							& $ 3.1536 \, 10^{-3} $   					& $ 0.20 $                  & $ 1.58 \, 10^{-3} $     	  \\ \hline
	    		dalleobtur 							& $ 3.1536 \, 10^{-3} $ 					& $ 0.20 $                  & $ 1.58 \, 10^{-3} $     	  \\ \hline
	    		drain 									& $ 94608 $   									& $ 0.30 $                  & $ 1 $     	  					\\ \hline
	    		conteneur1/conteneur2 		& $ 3.1536 \, 10^{-4} $   					& $ 0.12 $                  & $ 4.47 \, 10^{-4} $     	  \\ \hline
	    		dechet1/dechet2 					& $ 3.1536 \, 10^{-4} $   					& $ 0.30 $                  & $ 1.37 \, 10^{-3} $     	  \\ \hline	     		
\end{tabular}
\caption{[Test case 3] Data for flow and transport problems \cite{H17}.}
\label{A2Tab:ANDRAdata}  \vspace{-0.2cm}
\end{table}
For the spatial discretization of both flow and transport problems,  a non-uniform rectangular mesh is used as shown in Figure~\ref{A2Fig:ANDRAdomain} in blue, with $ 171 $ cells in the $ x-$direction and $ 158 $ cells in the $ y- $direction. The mesh size is $h \approx 0.42$m.  The Darcy flow problem~\eqref{A2ANDRADarcyflow} is solved by using the same mixed hybrid finite element method as presented in Section~\ref{sec:model}.  Numerical approximation of the hydraulic head is shown in Figure~\ref{A2Fig:ANDRADarcyDecom} (left).  For the transport problem, we decompose the domain into $ 6 $ rectangular subdomains in such a way that the black zone (\textit{terrain}) is separated from the rest and subdomain $ \Omega_{3} $ includes the \textit{dallerobtur}, \textit{voile}, \textit{radier} and a part of \textit{drain} zones (see Figure~\ref{A2Fig:ANDRADarcyDecom} (right)).
The transport is dominated by diffusion in subdomain $ \Omega_{3} $ (the maximum of the local P\'eclet number $\text{Pe}_{L} \approx  0.0032 $) and is dominated by advection (with $\text{Pe}_{L} \approx  2.75 $) in the other subdomains. 

\begin{figure}[!http]
\begin{minipage}[c]{0.95 \linewidth}
\centering
\includegraphics[scale=0.25]{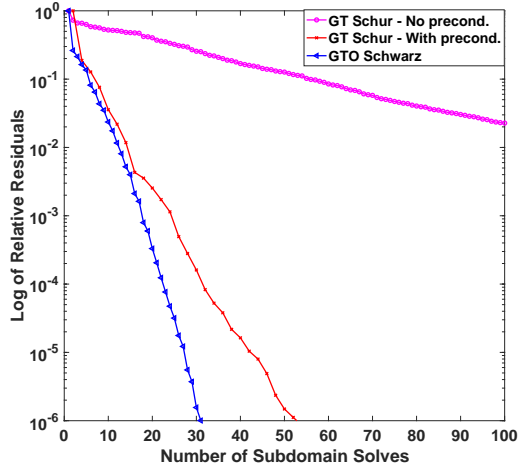} 
\end{minipage}  
\caption{[Test case 3] Relative residuals of GMRES for Method 1 (GT-Schur) with no preconditioner (magenta curves) and with the Neumann-Neumann preconditioner (red curves), and Method 2 (GTO-Schwarz) (blue curves).}
\label{A2Fig:ANDRArelres} \vspace{-0.2cm}
\end{figure}

For long time simulations,  we use time windows, i.e.  we split $ (0,T_{f}) $  into nonoverlapping smaller subintervals,  called time windows,  and then applies the DD methods in each time window successively in which the solution from the previous time window is used as the initial guess for the next time window.  For this test case, we use time windows of size $T=5$ years, and we will first analyze the convergence behavior as well as the accuracy in time of the multidomain solution with nonconforming grids for the first time window, $ (0,T) $.  The time steps are $ \Delta t_{3} =T/50, $ and $\Delta t_{i} = T/10 $, $ i \neq 3 $. 
%
%
The interface problem is solved iteratively using GMRES with a zero initial guess for both GT-Schur and GTO-Schwarz methods, the tolerance is set to be $10^{-6}$.
We show in Figure~\ref{A2Fig:ANDRArelres} the relative residuals for GT-Schur with or without preconditioning and GTO-Schwarz versus the number of subdomain solves.
We observe that the GTO-Schwarz method converges faster than the GT-Schur method with Neumann-Neumann preconditioner; without preconditioning,  the convergence of GT-Schur is very slow and it takes more than $450$ iterations to reach the same tolerance.  We remark that for this test case,  the advection is sufficient to make GTO-Schwarz faster than GTP-Schur with Neumann-Neumann preconditioning, but the advection term is small enough for the Neumann-Neumann preconditioner to be effective.

\begin{figure}[!http]
\centering
    \begin{subfigure}[b]{0.48\textwidth}
    \centering
        \includegraphics[width=\textwidth]{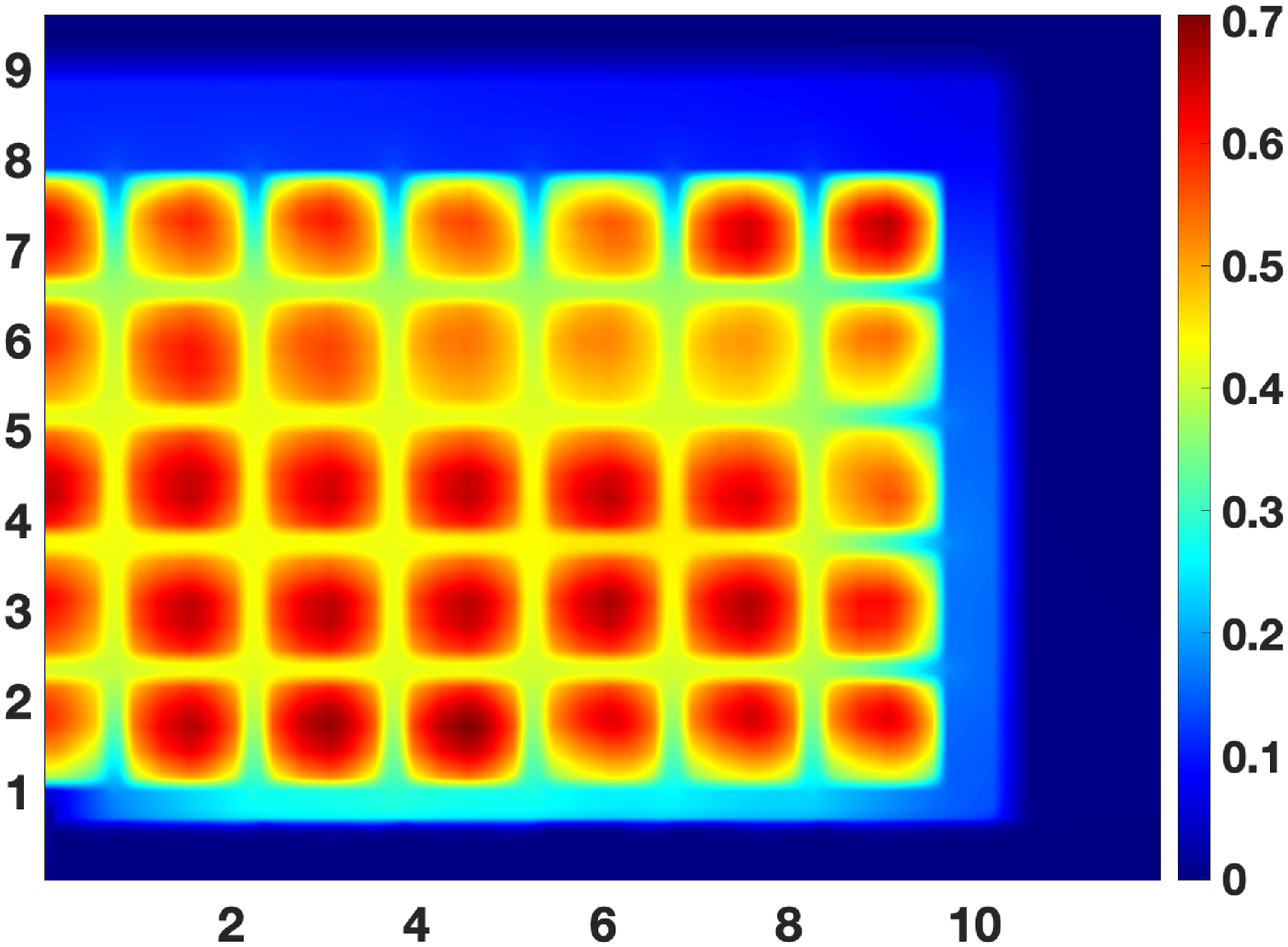}
    \end{subfigure}\hspace{-0.25cm}
    \begin{subfigure}[b]{0.48\textwidth}
     \centering
        \includegraphics[width=\textwidth]{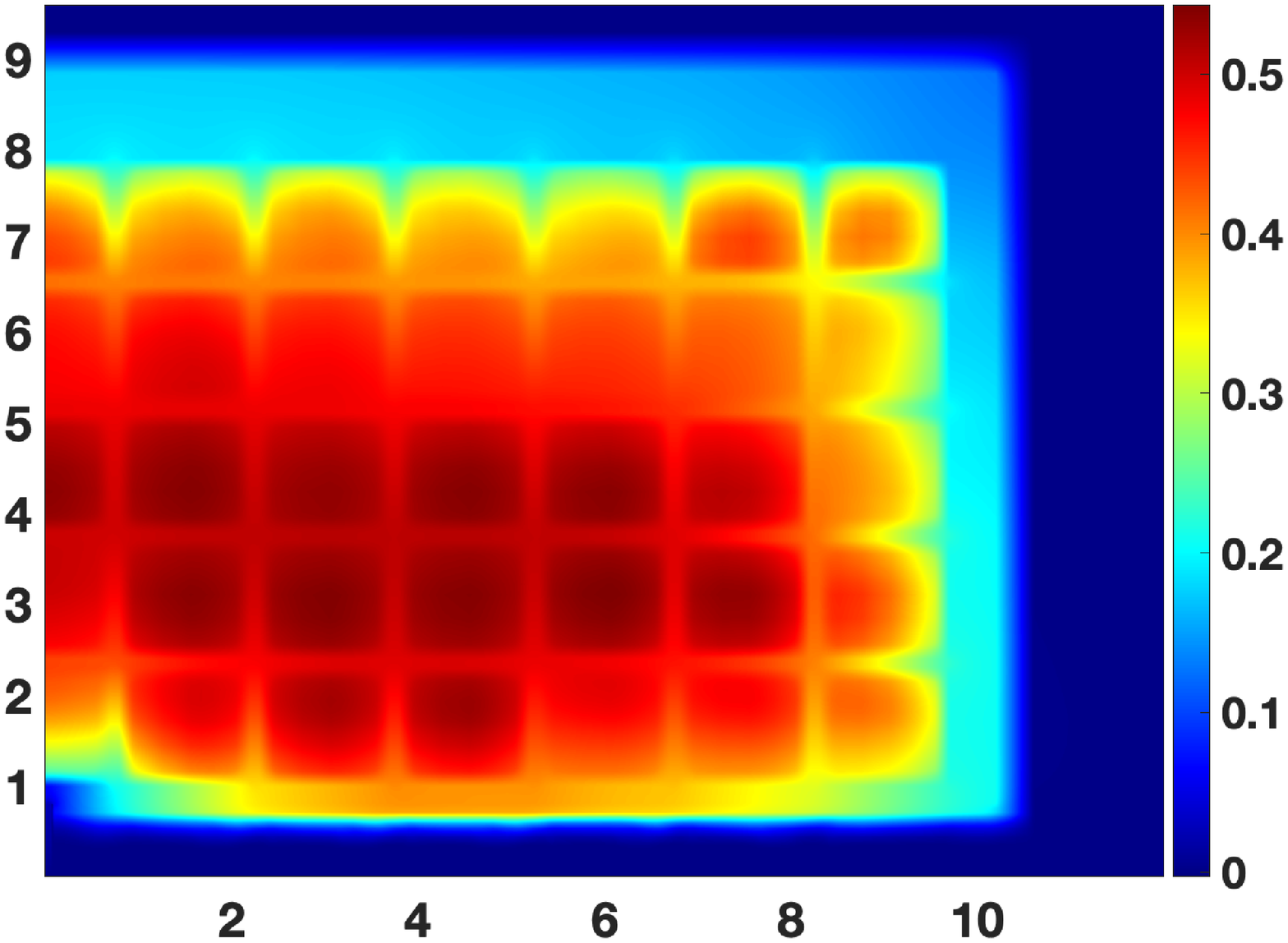}
    \end{subfigure}\hspace{-0.25cm}
    \begin{subfigure}[b]{0.48\textwidth}
     \centering
        \includegraphics[width=\textwidth]{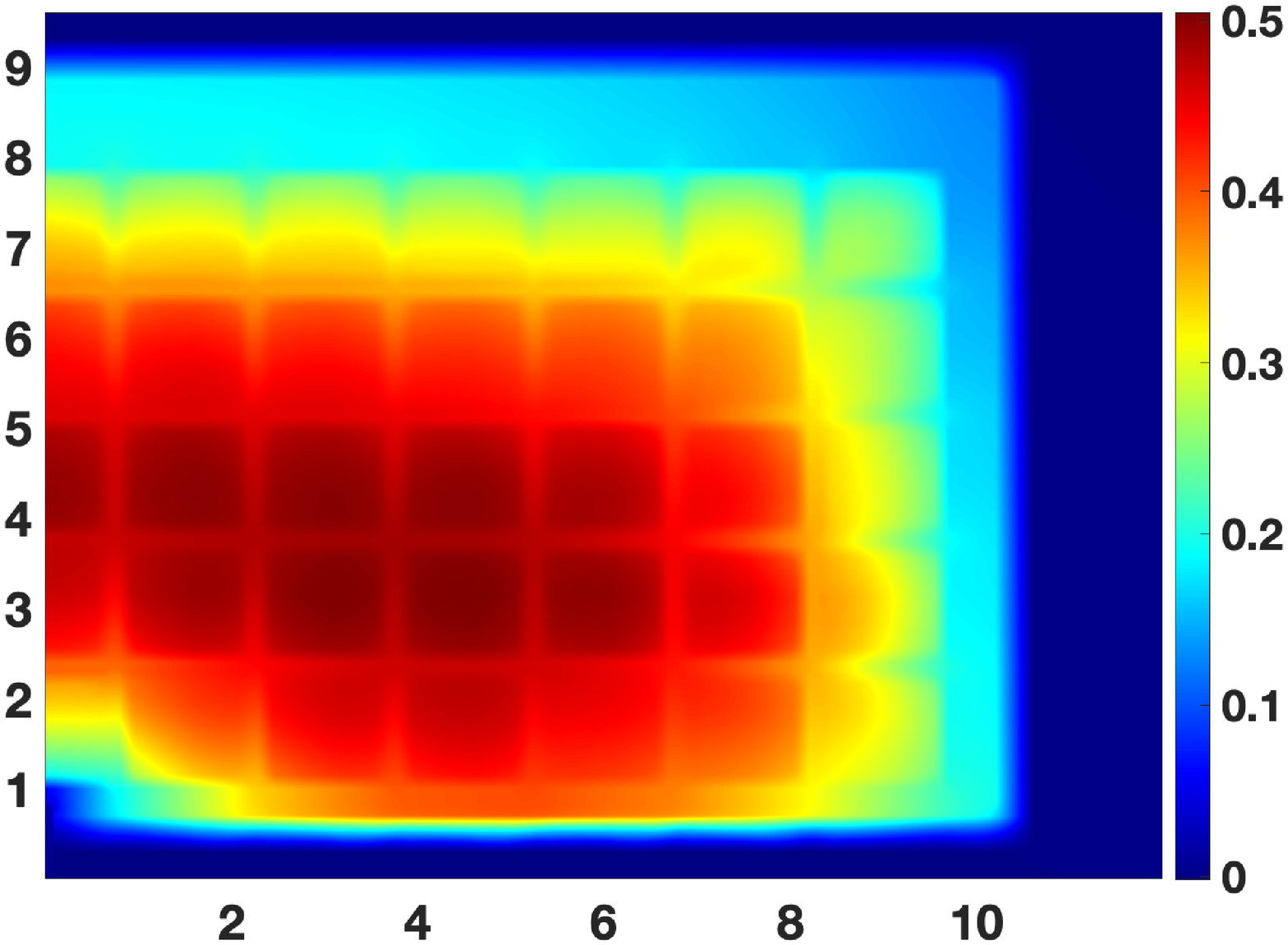}
    \end{subfigure} \hspace{-0.25cm}
    \begin{subfigure}[b]{0.48\textwidth}
     \centering
        \includegraphics[width=\textwidth]{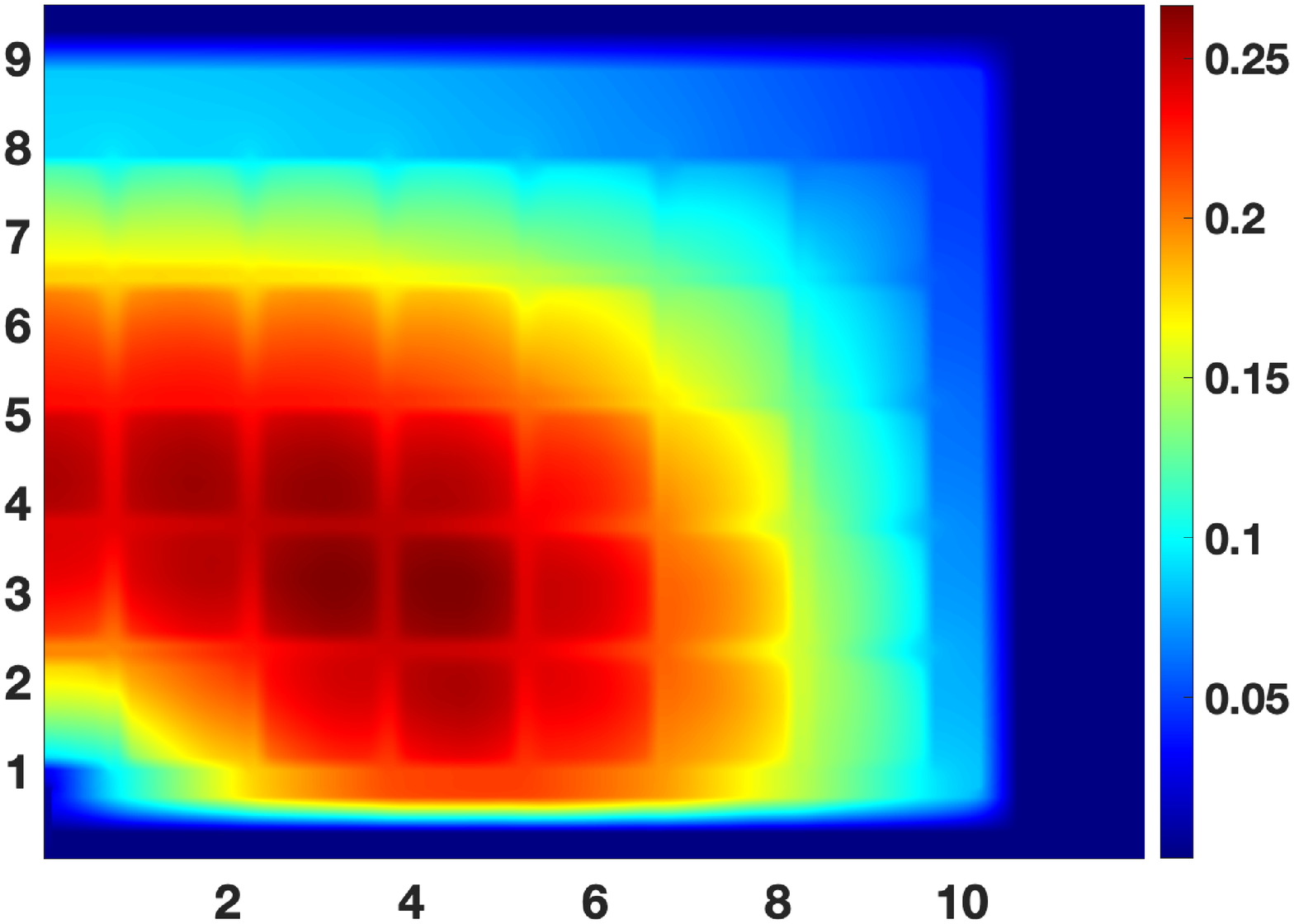}
    \end{subfigure} \vspace{-0.2cm}
    \caption{[Test case 3] Snapshots of the concentration after 20 years, 50 years, 100 years and 500 years respectively. The color bar for each plot is different. }\label{A2Fig:ANDRASolutionTW}
\end{figure}

%
%
Next, we run the GTP-Schur (with Neumann-Neumann preconditioner) and GTO-Schwarz methods for $ 100 $ time windows and stop the iterations in each time window when the relative residual is less than $ 10^{-3} $. 
The average number of iterations in each time window is approximately $ 8 $ (equivalent to $ 16 $ subdomain solves) for the GTP-Schur method and  is approximately $ 9 $ (equivalent to $ 9 $ subdomain solves) for the GTO-Schwarz method. 
Figure~\ref{A2Fig:ANDRASolutionTW} shows the concentration field after $20$~years,  $ 50 $~years, $100$~years and $ 500 $ years respectively (note that the color bar for each plot is different).  We see that the radionuclide escapes from the waste packages and slowly migrates into the surrounding area. Due to the specific design of the storage and under the effect of advection, the radionuclide tends to move toward the bottom right corner.

\section*{Conclusion} 
We have studied two global-in-time,  nonoverlapping DD methods for the linear advection-diffusion equation to model contaminant transport in heterogeneous porous media. The equation is discretized in space by a mixed hybrid method based on the lowest-order Raviart-Thomas finite element space with the flux variable consisting of both advective and diffusive flux.  Lagrange multipliers are introduced to enforce the continuity of the normal flux across the inter-element boundaries and are used for the discretization of the advective term.  The semi-discrete continuous-in-time problem is formulated as a space-time interface problem using either physical transmission conditions or Robin conditions,  which corresponds to the GTP-Schur with Neumann-Neumann preconditioner and the GTO-Schwarz method with optimized Robin parameters, respectively.  The developed methods are fully implicit in time and enable local time steps in the subdomains.  This work can be seen as a sequel to \cite{H13, H17} where similar methods were studied for the pure diffusion problem and the advection-diffusion problem with operator splitting, respectively.  Differently from \cite{H17},  here we do not treat advection and diffusive separately and no explicit time stepping is used.  We prove the convergence of the fully discrete OSWR algorithm with the upwind-mixed hybrid spatial discretization and backward Euler time-stepping method on nonconforming time grids.  Numerical results confirm the convergence and accuracy of the proposed methods with different time steps,  and their application to long-term simulations of transport of nuclear waste around a subsurface storage.  We observe that both methods handle well the case with large jumps in the coefficients, and their convergence is weakly dependent on the discretization parameters.  The GTP-Schur method works well and converges faster than without a preconditioner when the advection is moderate while the GTO-Schwarz method is insensitive to the advection and converges faster than the GTP-Schur method when there is sufficient advection. 
We are currently investigating the GTO-Schwarz method with second-order (Ventcell) transmission conditions \cite{H16Ventcell} and extension of the methods to the case of transport problems in fractured porous media as studied in \cite{H16} (for modeling the flow of a compressible fluid) where the fractures are treated as manifolds of one dimension less than the surrounding rock matrix.  


\bibliographystyle{elsarticle-num}

\end{document}